\documentclass{amsart}
\usepackage{righttag}
\usepackage{epsf}
\def\hepsffile{\leavevmode\epsffile}

\theoremstyle{plain}
\newtheorem{thm}{Theorem}[subsection]

\newtheorem{prop}[thm]{Proposition}

\theoremstyle{definition}
\newtheorem{defin}[thm]{Definition}

\newtheorem{emf}[thm]{}
\newtheorem{rem}[thm]{Remark}

\def\id{\protect\operatorname{id}}

\def\sign{\protect\operatorname{sign}}

\def\spin{\protect\operatorname{spin}}

\def\pr{\protect\operatorname{pr}}

\def\C{{\mathbb C}}
\def\Z{{\mathbb Z}}

\def\R{{\mathbb R}}
\def\N{{\mathbb N}}

\def\1{\hbox{\rm\rlap {1}\hskip.03in{\rom I}}}
\def\Bbbone{{\rm1\mathchoice{\kern-0.25em}{\kern-0.25em}
	{\kern-0.2em}{\kern-0.2em}I}}

\def\p{\partial}

\begin{document}
\hyphenation{Ca-m-po}

\title[Finite Order Invariants of Legendrian, Transverse, and Framed Knots]
{Finite Order Invariants of Legendrian, Transverse, and Framed Knots in Contact $3$-manifolds}
\author[V.~Tchernov]{Vladimir Tchernov}
\address{D-MATH, HG G66.4, ETH Zentrum, CH-8092 Z\"urich, Switzerland}
\email{chernov@math.ethz.ch}
\begin{abstract}
We show that for a big class of contact manifolds the groups of order $\leq
n$ invariants (with values in an arbitrary Abelian group) of Legendrian, of transverse,
and of framed knots are canonically isomorphic. 

On the other hand for an arbitrary cooriented contact structure on $S^1\times S^2$ 
with the nonzero Euler class of the contact bundle we construct examples of 
Legendrian homotopic Legendrian knots $K_1$ and $K_2$ such that they realize
isotopic framed knots but can be distinguished by finite order invariants of
Legendrian knots in $S^1\times S^2$. 
We construct similar examples for a big 
class of contact manifolds $M$ such that $M$ is a total space of a locally 
trivial $S^1$-fibration over a nonorientable surface.
We show that in some of these examples 
the complements of $K_1$ and of $K_2$ are overtwisted.

\end{abstract}
\maketitle

In this paper $\mathcal A$ is a not necessarily torsion free 
Abelian group, and $M$ is a not necessarily compact 
connected orientable three dimensional Riemannian manifold. 

\section{Introduction} 
A {\em contact structure\/} on a three-dimensional manifold $M$ is a smooth
field $\{ C_x\subset T_xM|x\in M\}$ of tangent two-dimensional planes,
locally defined as a kernel of a differential $1$-form $\alpha$ with non-vanishing
$\alpha\wedge d\alpha$. 
A contact structure is {\em cooriented\/} if the
two-dimensional planes defining the contact structure are continuously
cooriented (transversally oriented). A contact structure is 
{\em parallelized\/} if the two dimensional vector bundle 
$\{ C_x \}$ over $M$ is trivialized. 
A contact structure $C$ on a manifold $M$
is said to be {\em overtwisted\/} if there exists a $2$-disk $D$ embedded
into $M$ such that the boundary $\p D$ is tangent to $C$ while the disk $D$
is transverse to $C$ along $\p D$. Not overtwisted contact structures are
called {\em tight\/}.

A curve in $M$ is an immersion of $S^1$ into $M$. A curve in a contact manifold 
is called {\em Legendrian (resp. transverse)\/} if it is everywhere tangent 
(resp. transverse) to the planes of the contact structure. If the contact structure 
is cooriented, then every Legendrian curve has a natural framing. Similarly 
if the contact structure is parallelized, then every transverse curve 
has a natural framing. This means that every component of the space of 
Legendrian curves in a manifold $M$ with a cooriented contact structure 
is a subspace of a component of the space of framed curves.
Similarly every component of the space of transverse curves in a
manifold $M$ with a parallelized contact structure 
is a subspace of a component of the space of framed curves.
A Legendrian (resp. transverse) knot is a Legendrian (resp. transverse)
embedding of $S^1$ into a contact manifold.

In~\cite{FuchsTabachnikov} D.~Fuchs and S.~Tabachnikov proved that the 
quotients of the groups of $\C$-valued 
order $\leq n$ invariants by the groups of order $\leq (n-1)$ invariants 
of Legendrian, of transverse, and of framed knots in the standard
contact $\R^3$ are canonically isomorphic. An analogous
result was obtained by J.~W.~Hill~\cite{Hill} in the case of Legendrian and
of framed knots in the spherical cotangent bundle of $\R^2$ (with the
standard contact structure). 
To prove these facts they used the fact of the existence of the universal 
Vassiliev invariant of framed knots in these spaces (also known as
Kontsevich integral). Recently the existence of such an invariant was shown
by J.~E.~Andersen, J.~Mattes, and N.~Reshetikhin~\cite{AMR} for framed knots in
the spaces that are the product of $\R^1$ and an oriented compact surface with nonempty boundary.
However for most manifolds the existence of such a universal Vassiliev 
invariant is still unknown.
 
In this paper we use a different approach. It allows us to generalize
these isomorphisms to the case of finite order invariants
(with values in an arbitrary Abelian group) of Legendrian, of transverse, and of 
framed knots in a big class of contact manifolds. 
 
For a closed orientable manifold $M$ admitting a metric of a negative
sectional curvature we show that for an arbitrary cooriented contact structure 
on $M$ the groups of order $\leq n$ invariants with values in an arbitrary 
Abelian group of Legendrian and of framed knots (from the corresponding components
of the spaces of Legendrian and of framed curves) are canonically 
isomorphic.
For such manifolds with parallelized contact structure we show the same fact
for the groups of order $\leq n$ invariants of transverse and of framed knots.

The same fact is proven in the case of zero-homologous Legendrian,
transverse, and framed knots in an arbitrary 
contact manifold with the parallelized contact structure 

For a contact $3$-manifold $M$ (with the parallelized
contact structure) such that it is a total space of
a locally trivial $S^1$- or $\R^1$-fibration over a surface
we show that
the groups of order $\leq n$ invariants (with values in an arbitrary Abelian group) 
of Legendrian, of transverse, and of framed knots are canonically isomorphic,
provided that $M\neq S^1\times S^2$ and $M$ is not the quotient of 
$S^1\times S^2$ by the equivalence relation $e^{2\pi i t}\times x=e^{-2\pi i
t}\times (-x)$ (the last manifold $S^1$-fibers over $\R P^2$).

For an arbitrary cooriented contact structure on $S^1\times S^2$ with a
nonzero Euler class of the contact bundle we construct examples of 
Legendrian homotopic Legendrian knots $K_1$ and $K_2$ in $S^1\times S^2$ 
such that $K_1$ and $K_2$ realize isotopic framed knots but can be distinguished 
by order one invariants of Legendrian knots in $S^1\times S^2$. 

We also construct such examples for a big class of contact manifolds $M$ such
that $M$ is a total space of a locally trivial $S^1$-fibration over a
nonorientable surface. This means that in general
finite order invariants of Legendrian and of framed knots are different, 
and finite order invariants of Legendrian knots can distinguish
Legendrian knots that realize isotopic framed knots.

We show that in some of these examples the 
restrictions of the contact structure to the complements of 
$K_1$ and of $K_2$ are overtwisted. This seems to suggest that there is a 
nontrivial theory of Legendrian knots that realize isotopic
framed knots and have overtwisted complements.

\section{Definitions and known results}

\subsection{Definitions.}\label{definitions}
A {\em contact structure\/} on a three-dimensional manifold $M$ is a smooth
field $\{ C_x\subset T_xM|x\in M\}$ of tangent two-dimensional planes,
locally defined as a kernel of a differential $1$-form $\alpha$ with non-vanishing
$\alpha\wedge d\alpha$. A manifold with the contact structure possesses
canonical orientation determined by the volume form $\alpha\wedge\ d
\alpha$. The standard contact structure in $\R^3$ is given by the $1$-form
$\alpha=ydx-dz$.

A {\em contact element\/} on the manifold is a hyperplane
in the tangent space to the manifold at a point.
For a surface $F$ we denote by $ST^*F$ the space of all
cooriented (transversally oriented) contact elements of $F$. This space is a
spherical cotangent bundle of $F$. Its natural contact structure is a
distribution of tangent hyperplanes given by a condition
that a velocity vector of an incidence point of a contact element
belongs to the element.

A contact structure is {\em cooriented\/} if the
two-dimensional planes defining the contact structure are continuously
cooriented (transversally oriented). 
A contact structure is {\em oriented\/} if the two-dimensional planes defining the
contact structure are continuously oriented. Since every contact manifold
has a natural orientation we see that every cooriented contact structure is
naturally oriented and every oriented contact structure is naturally
cooriented.
A contact structure is {\em
parallelizable\/} ({\em parallelized\/}) if the two-dimensional vector bundle 
$\{ C_x \}$ over $M$ is trivializable (trivialized). Since every contact manifold
has a canonical orientation, one can see that every parallelized contact
structure is naturally cooriented. A contact structure $C$ on a manifold $M$
is said to be {\em overtwisted\/} if there exists a $2$-disk $D$ embedded
into $M$ such that the boundary $\p D$ is tangent to $C$ while the disk $D$
is transverse to $C$ along $\p D$. Not overtwisted contact structures are
called {\em tight\/}.

A {\em curve\/} in $M$ is an immersion of $S^1$ into $M$. 
(All curves have the natural orientation induced by the orientation of
$S^1$.) A {\em framed curve\/} in $M$ is a curve in $M$ equipped 
with the continuous unit normal vector field.

A {\em Legendrian curve\/} in a contact manifold $(M,C)$ is a curve
in $M$ that is everywhere tangent to $C$. If the contact structure on
$M$ is cooriented, then every Legendrian curve has a natural framing given 
by the unit normals to the planes of the contact structure that point in the
direction specified by the coorientation. 

For $(M,C)$ a contact manifold with a cooriented contact structure, 
we put $CM$ to be the total space of the fiberwise spherization of the
contact bundle, and we put $\pr:CM\rightarrow M$ to be the corresponding
locally trivial $S^1$-fibration.
The $h$-principle proved for the
Legendrian curves by M.~Gromov~\cite{Gromov} says that the space of Legendrian curves in
$(M,C)$ is weak homotopy equivalent to $\Omega CM$ the space of free loops
in $CM$. The equivalence is given by mapping a point of a Legendrian curve
to the point of $CM$ corresponding to the direction of the velocity vector
of the curve at this point.

To a Legendrian curve $K_l$ in the contact manifold with the parallelized 
contact structure one can associate an integer that is the
number of revolutions of the direction of the velocity vector of $K_l$ (with
respect to the chosen frames in $C$) under traversing $K_l$ according
to the orientation. This integer is called the {\em Maslov number\/} of
$K_l$. The $h$-principle 
implies that the set of Maslov numbers enumerates the
set of the connected components of the space of Legendrian 
curves in $\R^3$.

A {\em transverse\/} curve in a contact manifold $(M,C)$ is a curve in $M$
that is everywhere transverse to $C$. If the
contact structure on $M$ is parallelized, then a transverse curve has a
natural framing given by the unit normals corresponding to  
the projections on the tangent two-planes (orthogonal 
to the velocity vectors of the curve) of the first of the two coordinate 
vectors of the contact planes. A transverse curve in the contact manifold
with a cooriented contact structures is said to be {\em positive\/} if at
every point the velocity vector of the curve points into the coorienting
half-plane, and it is said to be {\em negative\/} otherwise.

A {\em knot ( framed knot)\/} in $M$ is an embedding (framed embedding) of $S^1$ into $M$.
In a similar way we define Legendrian and transverse knots in $M$.

A {\em singular (framed)\/} knot with $n$-double points is a curve (framed curve)
in $M$ whose only singularities are $n$ (transverse) double points.
An {\em isotopy\/} of a singular (framed) knot 
with $n$ double points is a path in the space of singular (framed) knots with
$n$ double points under which the preimages of the double points on $S^1$
change continuously.

An $\mathcal A$-valued framed (resp. Legendrian, resp. transverse) knot invariant is 
an $\mathcal A$-valued function on the set of the isotopy classes of framed (resp.
Legendrian, resp. transverse) knots.

A transverse double point $t$ of a singular knot can be resolved in two 
essentially different ways. We say that a resolution of a double point is
positive (resp. negative) if the tangent vector to the
first strand, the tangent vector to the second strand, and the vector from
the second strand to the first form the positive $3$-frame (this does not depend
on the order of the strands).
(If the singular knot is Legendrian (resp. transverse), 
then these resolution can be made in the
category of Legendrian (resp. transverse) knots.)
A singular framed (resp. Legendrian, resp. transverse) knot $K$ with $(n+1)$ 
transverse double points
admits $2^{(n+1)}$ possible resolutions of the double points. A sign of the resolution 
is put to be $+$ if the number of negatively resolved double points is even, and
it is put to be $-$ otherwise. 
Let $x$ be an $\mathcal A$-valued invariant of framed (resp. Legendrian,
resp. transverse) knots. The invariant $x$ is said to be of {\em finite
order\/} (or {\em Vassiliev invariant\/}) if there exists a positive 
integer $(n+1)$ such that for any singular knot $K_s$ with $(n+1)$
transverse double points the sum (with appropriate signs) of the values of $x$ on the nonsingular
knots obtained by the $2^{n+1}$ resolutions of the double points is zero. 
An invariant is said to be of order not greater than $n$ (of order $\leq n$) if $n$
can be chosen as integer in the above definition. The group of $\mathcal
A$-valued finite order invariants has an increasing filtration by the
subgroups of the invariants of order $\leq n$.

\begin{emf}\label{description}{\em Description of Legendrian and of transverse
knots in $\R^3$.\/} 
The contact Darboux theorem says that every contact manifold locally looks
like $\R^3$ with the standard contact structure $\alpha=ydx-dz$. 
Transverse and Legendrian knots are conveniently presented by the projections 
into the plane $(x,z)$. Identify a  point $(x, y, z)\in \R^3$ with the
point $(x,z)\in \R^2$ furnished with a fixed direction of an unoriented 
straight line through $(x,z)$ with the slope $y$. Then the curve in $\R^3$
is a one parameter family of points with directions in $\R^2$. 

A curve in
$\R^3$ is transverse if and only if the corresponding curve in $\R^2$ is
never tangent to the chosen directions along itself.

While a generic regular curve has a regular projection onto the
$(x,z)$-plane, the projection of a generic Legendrian curve onto the
$(x,z)$-plane has isolated critical points (since all the planes of the
contact structure are parallel to the $y$-axis). Hence the projection of a
generic Legendrian curve may have cusps. A curve in $\R^3$ is Legendrian if
and only if the corresponding planar curve with cusps 
is everywhere tangent to the field of directions, in particular this field
is determined by the curve with cusps.
\end{emf}

\subsection{Generalized Theorem of D.~Fuchs and S.~Tabachnikov for
Legendrian knots.}

Let $(M,C)$ be a contact manifold with a cooriented contact structure.
Let $\mathcal L$ be a component of the space of Legendrian curves in
$(M,C)$. 
Since a Legendrian curve in $(M,C)$ has a natural framing, we see that
$\mathcal L$ is a subspace of a connected component $\mathcal F$ of the
space of framed curves in $M$.

In~\cite{FuchsTabachnikov} D.~Fuchs and S.~Tabachnikov 
proved that $\C$-valued finite order
invariants of Legendrian knots in the standard contact $\R^3$ can not distinguish 
Legendrian knots that belong to the same component $\mathcal L$ of the space of 
Legendrian curves in $\R^3$ and realize isotopic framed knots. 

\begin{emf}\label{conditions} 
As it was later observed by D.~Fuchs and S.~Tabachnikov~\cite{FuchsTabachnikovprivate},
the proof of the Theorem can be generalized to a large class of contact
manifolds with cooriented contact structures. Moreover it can be generalized
to the case of finite order invariants with values in an arbitrary Abelian
group $\mathcal A$. In order to state the generalized theorem of D.~Fuchs
and S.~Tabachnikov we introduce the following definitions.

Let $(M,C)$ be a contact manifold with a cooriented contact structure, let
$\mathcal L$ be a component of the space of Legendrian curves in $(M, C)$,
and let $\mathcal F$ be the corresponding component of the space of framed curves
in $M$. 
We say that $\mathcal F$ {\em satisfies condition\/} $\textrm{I}$ if
for any $0\neq n\in\N$ every framed knot $K\in\mathcal F$ is not isotopic to
the framed knot $K_n$ obtained by adding $2n$ positive extra 
twists to the framing of $K$. 
We say that $\mathcal F$ {\em satisfies condition\/} $\textrm{II}$ 
if it contains infinitely many components of the space of Legendrian curves
(see Proposition~\ref{interpretationconditionII} 
%and Remark~\ref{remarkinterpretation} 
for the homological interpretation of condition \textrm{II}).

\end{emf} 

Below we formulate the generalized theorem of D.~Fuchs and S.~Tabachnikov
for Legendrian knots.

\begin{thm}\label{thmTF}
Let $(M,C)$ be a contact manifold
with a cooriented contact structure. Let $\mathcal L$ be a component of
the space of Legendrian curves in $(M,C)$, and let $\mathcal F$ be the
corresponding component of the space of framed curves in $M$.
Let $x$ be an $\mathcal A$-valued order $\leq n$
invariant of Legendrian knots in $(M, C)$, and let $K_1,K_2\in \mathcal L$ 
be Legendrian knots that represent isotopic framed knots.

Then $x(K_1)=x(K_2)$, provided that $\mathcal F$ satisfies conditions
\textrm{I} and \textrm{II}.
\end{thm}

The proof of this Theorem is a more or less straightforward generalization 
of the proof of the corresponding Theorem for $\mathcal A=\C$ and $M=\R^3$
proved by D.~Fuchs and S.~Tabachnikov~\cite{FuchsTabachnikov}.

\subsection{Generalized Theorem of D.~Fuchs and S.~Tabachnikov for
transverse knots}

Let $(M,C)$ be a contact manifold with a parallelized contact structure.
Let $\mathcal T$ be a component of the space of transverse curves in
$(M,C)$. 
Since a transverse curve in $(M,C)$ has a natural framing, we see that
$\mathcal T$ is a subspace of a connected component $\mathcal F$ of the
space of framed curves in $M$.

In~\cite{FuchsTabachnikov} D.~Fuchs and S.~Tabachnikov 
proved that $\C$-valued finite order
invariants of transverse knots in the standard contact $\R^3$ can not distinguish 
transverse knots that belong to the same component $\mathcal T$ of the space of 
transverse curves in $\R^3$ and realize isotopic framed knots. 

\begin{emf}
As it was later observed by D.~Fuchs and S.~Tabachnikov~\cite{FuchsTabachnikovprivate},
the proof of the Theorem can be generalized to a large class of contact
manifolds with parallelized contact structures. Moreover it can be generalized
to the case of finite order invariants with values in an arbitrary Abelian
group $\mathcal A$. In order to state the generalized theorem of D.~Fuchs
and S.~Tabachnikov we introduce the following definition.

Let $(M,C)$ be a contact manifold with a parallelized contact structure, let
$\mathcal T$ be a component of the space of transverse curves in $(M, C)$,
and let $\mathcal F$ be the corresponding component of the space of framed curves
in $M$. 
We say that $\mathcal F$ {\em satisfies condition\/} $\textrm{I}$ if
for any $0\neq n\in\N$ every framed knot $K\in\mathcal F$ is not isotopic to
the framed knot $K_n$ obtained by adding $2n$ positive extra 
twists to the framing of $K$. 
\end{emf} 

Below we formulate the generalized theorem of D.~Fuchs and S.~Tabachnikov
for transverse knots. 

\begin{thm}\label{thmTFtransverse} 
Let $(M,C)$ be a contact manifold with a parallelized contact structure. Let 
$\mathcal T$ be a component of the space of transverse curves in $(M,C)$, and 
let $\mathcal F$ be the corresponding component of the space of framed curves 
in $M$. Let $x$ be an $\mathcal A$-valued order $\leq n$
invariant of transverse knots in $(M, C)$. Let $K_1,K_2\in \mathcal T$ 
be transverse knots that represent isotopic framed knots.

Then $x(K_1)=x(K_2)$, provided that $\mathcal F$ satisfies condition
\textrm{I}.
\end{thm}

The proof of this Theorem is a more or less straightforward generalization 
of the proof of the corresponding Theorem for $\mathcal A=\C$ and $M=\R^3$
proved by D.~Fuchs and S.~Tabachnikov~\cite{FuchsTabachnikov}.

\section{Isomorphism between the groups of order $\leq n$ invariants of
Legendrian and of framed knots}

Let $(M,C)$ be a contact manifold with a cooriented contact structure.
Let $\mathcal L$ be a connected component of
the space of Legendrian curves in $M$, and let $\mathcal F$ be the
connected
component of the space of framed curves that contains $\mathcal L$. 
Let $V_n^{\mathcal L}$ (resp. $W_n^{\mathcal F}$) 
be the groups of $\mathcal A$-valued order $\leq n$ invariants of Legendrian
(resp. framed) knots from $\mathcal L$ (resp.
from $\mathcal F$). Clearly every invariant $y\in W_n^{\mathcal F}$ restricted to
the category of Legendrian knots in $\mathcal L$ is an element $\phi(y)\in
V_n^{\mathcal L}$. This gives a homomorphism $\phi:W_n^{\mathcal
F}\rightarrow V_n^{\mathcal L}$.

\begin{thm}\label{isomorphism} Let $(M, C)$ be a contact manifold with a 
cooriented contact structure. Let $\mathcal L$ be a connected component of
the space of Legendrian curves in $M$, and let $\mathcal F$ be the
connected
component of the space of framed curves that contains $\mathcal L$.
Then the following two statements {\bf a:} and {\bf b:} are equivalent.
\begin{description}
\item[a] For any $x\in V_n^{\mathcal L}$ and any knots $K_1, K_2\in
\mathcal L$ that represent isotopic framed knots 
$x(K_1)=x(K_2)$.
\item[b] 
For $x\in V_n^{\mathcal L}$ there exists $\psi (x)\in W_n^{\mathcal F}$ such 
that $\phi(\psi(x))=x$. Such $\psi (x)$ is unique and it defines a
canonical isomorphism $\psi:V_n^{\mathcal L}\rightarrow W_n^{\mathcal F}$.
\end{description}
\end{thm} 

For the proof of Theorem~\ref{isomorphism} see
Subsection~\ref{proofisomorphism}.

The following Theorem is an immediate consequence of Theorems~\ref{thmTF}
and~\ref{isomorphism}.

\begin{thm}\label{isomorphismobtained}
Let $(M, C)$ be a contact manifold with a
cooriented contact structure, and let $\mathcal L$ be a connected component of
the space of Legendrian curves in $M$. Let $\mathcal F$ be the
connected component of the space of framed curves that contains $\mathcal L$.
Let $V_n^{\mathcal L}$ (resp. $W_n^{\mathcal F}$) 
be the groups of $\mathcal A$-valued order $\leq n$ invariants of Legendrian
(resp. framed) knots from $\mathcal L$ (resp. 
from $\mathcal F$). 
Then the groups $V_n^{\mathcal L}$ and $W_n^{\mathcal F}$ are canonically
isomorphic, provided that $\mathcal F$ satisfies conditions \textrm{I} and
\textrm{II}.
\end{thm}

\begin{rem}
Let $\mathcal F$ be a connected component of the space of framed curves in
$M$ that satisfies condition \textrm{I}. 
Theorem~\ref{isomorphismobtained} implies for any $n\in\N$ 
the group of $\mathcal A$-valued order 
$\leq n$ invariants of Legendrian knots from a connected component 
of Legendrian curves curves contained in $\mathcal F$ does not depend on the
choice of a cooriented contact structure, provided that for this choice
$\mathcal F$ satisfies condition \textrm{II}. And hence this group can not
be used to distinguish such cooriented contact structures.

The $h$-principle proved for the Legendrian curves by
M.~Gromov~\cite{Gromov} implies that $\mathcal F$ satisfies condition
\textrm{II} provided that the contact structure is parallelizable.
Clearly if $\mathcal F$ consists of zero homologous knots, 
then it satisfies condition \textrm{I} and we get that for an arbitrary 
contact manifold $(M,C)$ (with the parallelizable contact structure) and 
for an arbitrary component $\mathcal L$
consisting of zero homologous Legendrian curves the groups $V_n^{\mathcal
L}$ and $W_n^{\mathcal F}$ are canonically   
isomorphic. 

Since $H_1(\R^3)=0$ and the standard contact structure on $\R^3$ is
parallelizable, we get this isomorphism for the components
of the space of Legendrian and of framed curves in the standard contact $\R^3$.
In~\cite{FuchsTabachnikov} D.~Fuchs and S.~Tabachnikov showed that for the
standard contact $\R^3$ and for $\mathcal A=\C$ the quotient groups 
$V_n^{\mathcal L}/V_{n-1}^{\mathcal L}$ and $W_n^{\mathcal
F}/W_{n-1}^{\mathcal F}$ are canonically isomorphic. Their proof was based
on the fact that for the $\C$-valued Vassiliev invariants of framed knots in
$\R^3$ there exists the universal Vassiliev invariant 
constructed by T.~Q.~T.~Le and J.~Murakami~\cite{LeMurakami}. 
(For unframed knots in $\R^3$ the existence of the universal Vassiliev 
invariant is the classical result of M.~Kontsevich~\cite{Kontsevich}.)
Our results generalize the results of Fuchs and Tabachnikov.
\end{rem} 

\begin{thm}\label{fibration} 
Let $M$ be an orientable manifold that is a total space of a locally
trivial $S^1$-fibration over a (not necessarily closed or orientable) 
surface $F$. Assume moreover that $M\neq S^1\times S^2$ and $M$ is not the
quotient of $S^1\times S^2$ by the equivalence relation that $(e^{2\pi i
t}\times x)=(e^{-2 \pi i t}\times (-x))$ (such $M$ fibers over $\R P^2$). 
Then all the components of the space of framed curves in $M$ satisfy
condition \textrm{I}. 
\end{thm}

For the Proof of Theorem~\ref{fibration} see
Subsection~\ref{Prooffibration}.

\begin{emf}\label{corollariesfibration}{\em Corollaries of 
Theorem~\ref{fibration}.\/} 
Let $M$ be a manifold satisfying conditions of
Theorem~\ref{fibration}, let $C$ be a cooriented contact structure on $M$,
and let $\mathcal F$ be a connected component of the space of framed curves 
in $(M,C)$ that satisfies condition \textrm{II}. Then for every connected component $\mathcal L\subset
\mathcal F$ of the space of Legendrian curves there is a natural isomorphism
between the groups $V_n^{\mathcal L}$ and $W_n^{\mathcal F}$ 
of $\mathcal A$-valued order $\leq n$ invariants. 

Thus if the contact structure on a
manifold (satisfying conditions of Theorem~\ref{fibration}) is
parallelizable, then we get the isomorphisms of the groups of $\mathcal A$-valued 
order $\leq n$ invariants of Legendrian and of framed knots (from the
corresponding components of the two spaces).

A very interesting class of contact manifolds satisfying the conditions of
Theorem~\ref{fibration} is formed by the spherical
cotangent bundles $ST^*F$ of surfaces $F$ with the natural contact
structure on $ST^*F$ (see~\ref{definitions}). The theory of the invariants 
of Legendrian knots in $ST^*F$ is often referred to as the theory of 
V.~I.~Arnold's~\cite{Arnoldsplit} $J^+$-type invariants of fronts on a surface $F$. 
The natural contact structure $C$ on the spherical cotangent bundle $ST^*F$ of
a surface $F$ is cooriented. (The coorientation is induced from the
coorientation of the contact elements of $F$.) 
One can verify that for orientable $F$ 
the standard contact structure on $ST^*F$ is parallelizable, and hence all
the components of the space of framed curves satisfy condition~\textrm{II}.
If $F$ is not orientable, then the standard cooriented contact structure on
$ST^*F$ is not parallelizable, but one can still verify 
(cf. Proposition 8.2.4~\cite{Tchernov}) that every connected component of the
space of framed curves satisfies condition~\textrm{II}. 
Hence for an arbitrary surface $F$ we obtain the canonical isomorphism of   
the groups of $\mathcal A$-valued order $\leq n$ invariants of Legendrian 
and of framed knots (from the corresponding components of the spaces of 
Legendrian and of framed knots in $ST^*F$ with the standard contact structure). 
Or equivalently we get that the groups of $\mathcal A$-valued order 
$\leq n$ $J^+$-type invariants of fronts on $F$
and the groups of $\mathcal A$-valued order $\leq n$ invariants of framed
knots in $ST^*F$ (from the corresponding components of the two spaces) are canonically
isomorphic.

In~\cite{Hill} J.~W.~Hill showed that for the
standard contact $ST^*\R^2$ and for $\mathcal A=\C$ the quotient groups 
$V_n^{\mathcal L}/V_{n-1}^{\mathcal L}$ and $W_n^{\mathcal
F}/W_{n-1}^{\mathcal F}$ are canonically isomorphic. His proof was based
on the fact that for the $\C$-valued Vassiliev invariants of framed knots in
$ST^*\R^2$ there exists the universal Vassiliev invariant 
constructed by V.~Goryunov~\cite{Goryunov}.
(For unframed knots in $\R^3$ the existence of the universal Vassiliev 
invariant is the classical result of M.~Kontsevich~\cite{Kontsevich}.)
Our results generalize the results of J.~W. Hill (even in the case of
$M=ST^*\R^2)$.
\end{emf}

Before we proceed further we clarify what it means that a connected
component $\mathcal F$ of the space of framed curves in a contact manifold
$(M, C)$ with a cooriented contact structure satisfies condition \textrm{II}.

\begin{prop}\label{interpretationconditionII}
Let $(M,C)$ be a 
%closed 
contact manifold with a cooriented contact
structure, and let $\mathcal F$ be a component of the space of framed curves
in $M$. Then $\mathcal F$ does not satisfy condition \textrm{II} if and only if 
there exists $\alpha\in H_2(M, \Z)$ such that the value of the Euler class of the 
contact bundle on $\alpha$ is nonzero and $\alpha$ is realizable by a mapping 
$\mu:T^2\rightarrow M$ of the two torus with the property that $\mu$ 
of the meridian of $T^2$ is a loop free homotopic to loops realized 
by curves from $\mathcal F$.
\end{prop}

For the Proof of Proposition~\ref{interpretationconditionII} see
Subsection~\ref{ProofinterpretationconditionII}.

%\begin{rem}\label{remarkinterpretation}
%Using the proof of Proposition~\ref{interpretationconditionII}
%one can easily verify that if $(M,C)$ is a not closed contact manifold with
%a cooriented contact structure and $\mathcal F$ is a connected component 
%of the space of framed curves. Then $\mathcal F$ does not 
%satisfy condition \textrm{II} if and only if there exists a mapping
%$\mu:T^2\rightarrow M$ with the property that $\mu$ of the meridian of
%$T^2$ is free homotopic to loops realized by curves from $\mathcal F$ and
%such that $e_{\mu}\neq 0\in \Z=H^2(T^2, Z)$, where 
%$e_{\mu}$ is the Euler class of the oriented two-dimensional bundle over $T^2$ 
%obtained by a pullback (via $\mu$) of the contact bundle of $M$.
%\end{rem}

\begin{thm}\label{atoroidal}
Let $(M,C)$ be a contact manifold (with a cooriented contact
structure) such that $\pi_2(M)=0$, every $\alpha\neq 1\in\pi_1(M)$ is of 
infinite order, and for an arbitrary mapping
$\mu:T^2\rightarrow M$ of the two-torus 
the homomorphism $\mu_*:\pi_1(T^2)\rightarrow \pi_1(M)$ is not injective.
Then all the components of the space of framed curves in $M$ satisfy
conditions \textrm{I} and \textrm{II}, and hence the groups of $\mathcal
A$-valued order $\leq n$ invariants of Legendrian and of framed knots (from the
corresponding components) are canonically isomorphic.
\end{thm} 

For the Proof of Theorem~\ref{atoroidal} see
Subsection~\ref{Proofatoridal}.

\begin{rem} A well known Theorem by A.~Preissman (see~\cite{Docarmo} pp. 258-265) 
says that every nontrivial commutative subgroup of the fundamental group of a 
closed three-dimensional manifold $M$ of negative sectional curvature is infinite 
cyclic. It is also known that the universal covering of such $M$ is diffeomorphic 
to $\R^3$, and hence $\pi_2(M)=0$. Thus every closed manifold $M$ 
admitting a metric of a negative  sectional 
curvature satisfies all the conditions of Theorem~\ref{atoroidal} and for an arbitrary 
cooriented contact structure on such $M$ we obtain the isomorphism of the groups 
of $\mathcal A$-valued order $\leq n$ invariants of Legendrian and of framed
knots from the corresponding components.

Theorem~\ref{atoroidal} implies (see Proposition~\ref{Preissman})
that for $(M, C)$ a contact manifold with a cooriented contact structure 
such that $M$ admits a structure of an $\R^1$-fibration over 
$F\neq S^2, T^2, \R P^2, K$ (the Klein bottle) we also get the isomorphism of 
the groups of the $\mathcal A$-valued order $\leq n$ invariants of Legendrian 
and of framed knots from the corresponding components of the two spaces.
(From the work of T.~Fiedler~\cite{Fiedler} one can easily see that if $M$ 
is an orientable manifold $\R^1$-fibered over a (not necessarily orientable) 
surface $F$, then all the components of the space of framed curves satisfy 
condition \textrm{I}. Hence for an arbitrary parallelizable contact structure on 
such $M$ the groups of $\mathcal A$-valued order $\leq n$ invariants of Legendrian 
and of framed knots (from the corresponding components) are canonically
isomorphic.) 
\end{rem}

\subsection{Examples of Legendrian knots that are distinguishable by finite
order invariants.}
In this subsection we construct a big class of contact manifolds $(M, C)$ and
connected components $\mathcal L$ of Legendrian curves in $(M, C)$ 
for which there exist Legendrian knots $K_1, K_2\in \mathcal L$ such that 
$K_1$ and $K_2$ realize isotopic framed knots but can be distinguished by 
finite order invariants of Legendrian knots.

The Theorem of R.~Lutz~\cite{Lutz} says that for an arbitrary
orientable three-manifold $M$ every homotopy class of the distribution of
tangent two-planes to $M$ contains a contact structure.
(The Theorem of Ya. Eliashberg~\cite{Eliashberg} says even more that every homotopy 
class of the distribution of tangent two-planes to $M$ contains a positive
overtwisted contact structure.)

However in our constructions we will use only the Euler classes of the
contact bundles. For this reason we start with the following Proposition.

\begin{prop}\label{existcontact} 
\begin{description}
\item[1] Let $(M, C)$ be a 
%closed 
contact manifold with a cooriented contact
structure, and let $e\in H^2(M, \Z)$ be the Euler class of the contact bundle.
Then $e=2\alpha$, for some $\alpha\in H^2(M, \Z)$.
\item[2] Let $M$ be an 
%closed 
oriented manifold, and $e\in H^2(M, \Z)$ be
such that $e=2\alpha\in H^2(M, \Z)$ for some $\alpha\in H^2(M)$. 
Then there exists a cooriented contact structure $C$ on $M$ such that the 
the Euler class of the contact bundle is $e$.
\end{description}
\end{prop}

For the Proof of Proposition~\ref{existcontact} see
Subsubsection~\ref{Proofexistcontact}.

\begin{figure}[htbp]
 \begin{center}
  \epsfxsize 8cm
  \hepsffile{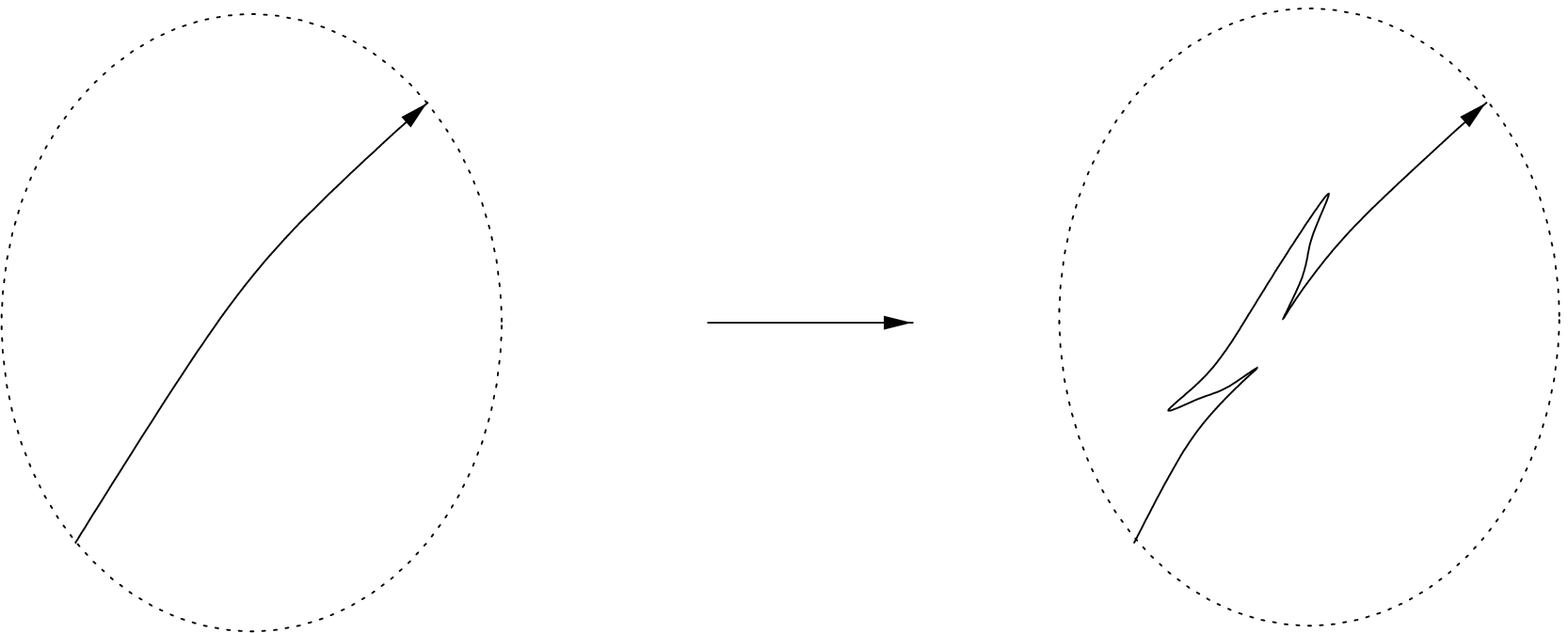}
 \end{center}
\caption{}\label{change.fig}
\end{figure}

Let $C$ be a cooriented contact structure on $M=S^2\times S^1$ such that
the Euler class of the contact bundle is nonzero.
Let $K$ be a knot in $S^2\times S^1$ that crosses exactly once one of the
spheres $S^2\times {t}$. The Theorem of Chow~\cite{Chow} and
Rashevskii~\cite{Rashevskii} says
that there exists a Legendrian knot $K_0$ that is $C^0$-small 
isotopic to $K$ as an unframed knot. Let $K_1$ be a Legendrian knot that
is the same as $K_1$ except of a small piece located in a chart
contactomorphic to the standard contact $\R^3$ where it is changed as it is
shown in Figure~\ref{change.fig} (see~\ref{description}).

\begin{thm}\label{example1} 
\begin{description}
\item[a] Legendrian knots $K_0$ and $K_1$ belong to the same component of
the space of Legendrian curves and realize isotopic framed knots.
\item[b] There exists a $\Z$-valued order one invariant $I$ of Legendrian
knots, such that $I(K_0)\neq I(K_1)$.
\end{description}
\end{thm}

For the Proof of Theorem~\ref{example1} see Subsection~\ref{Proofexample1}.

\begin{emf}\label{overtwisted}{\bf The contact structure and the knots $K_0$
and $K_1$ in Theorem~\ref{example1} can be
chosen so that the complements of $K_0$ and of $K_1$ are overtwisted.}
Let $\Delta$ be an embedded into $M$ disk centered at a point $p\in M$. The
Theorem of Eliashberg~\cite{Eliashberg} says that every homotopy class of
the distribution of two-planes contains an overtwisted contact structure
that has $\Delta$ as the standard overtwisted disk. In the example of
Theorem~\ref{example1} we can start with an unframed knot $K$ that is far away
from $\Delta$. Then since both $K_0$ and $K_1$ were constructed using a
$C^0$-small approximation of $K$ we can assume that they are
also far away from $\Delta$ and thus the restrictions of the contact
structure to the complements of $K_0$ and of $K_1$ are overtwisted. 
\end{emf}

\begin{rem}

Let $K_i$, $i\in\N$, be the knot that is the same as $K_0$ 
everywhere except of a small piece located in a chart contactomorphic 
to the standard contact $\R^3$ where it is changed in the way described 
by the addition of $i$ zigzags shown in Figure~\ref{change.fig}. 
The Proof of Theorem~\ref{example1} implies that all $K_i$'s are Legendrian 
homotopic knots that realize isotopic framed knots, but for all $i_1\neq i_2$ 
Legendrian knots $K_{i_1}$ and $K_{i_2}$ are not isotopic and there exists 
a $\Z$-valued order one invariant $I$ of Legendrian knots such that 
$I(K_{i_1})=I(K_{i_2})+(i_2-i_1)$. (And hence this $I$ distinguishes all the 
$K_i$'s.)
\end{rem}

\begin{emf}\label{OtherExamples}
{\em Below we describe another big family of examples when finite
order invariants distinguish Legendrian knots that realize isotopic framed
knots.\/} 

Let $F$ be a nonorientable surface that can be decomposed as a sum of
the Klein bottle $K$ and a surface $F'\neq S^2$.
Let $M$ be an orientable manifold that admits a structure of a locally
trivial $S^1$-fibration $p:M\rightarrow F$. (For example one can take $M$ to be $STF$
the spherical tangent bundle of $F$.)

\begin{figure}[htbp]
 \begin{center}
  \epsfxsize 6cm
  \hepsffile{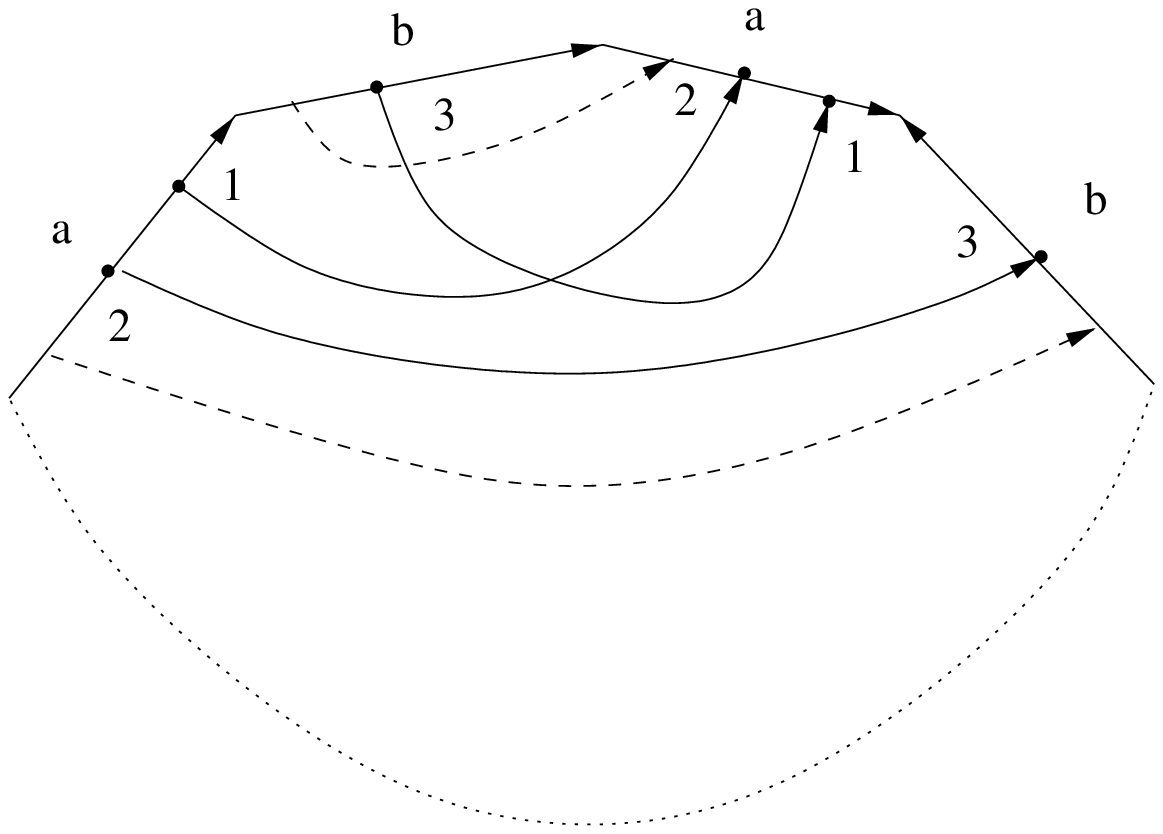}
 \end{center}
\caption{}\label{example2.fig}
\end{figure}

Consider an $S^1$-fibration $\xi:M'\rightarrow S^1$ induced from $p$ by the
mapping of $S^1$ to the solid loop in the surface shown in Figure~\ref{example2.fig}. 
(In this Figure the enumeration of the end points of the arcs indicates which pairs 
of points should be identified to obtain the loop.) 
Since the solid loop is an orientation preserving loop in $F$ we get that $M'=T^2$. 
Consider the natural map $\mu:M'=T^2\rightarrow M$ (the one from the definition 
of the induced fibration). Since a homology class in $H_1(M, \Z)$ projecting
to the dashed loop in Figure~\ref{example2.fig} has a nonzero intersection 
with the class realized by $\mu(T^2)$, we get that the class realized by
$\mu(T^2)$ is nonzero in $H_2(M, \Z)$. 

Let $C$ be a cooriented contact structure on $M$ such that the value of the
Euler class $e\in H^2(M, \Z)$ of the contact bundle on the homology class
realized by $\mu(T^2)$ is nonzero. And hence is equal to $2r$ for some
$0\neq r\in\Z$, see Proposition~\ref{existcontact}.

Let $K$ be an arbitrary Legendrian knot whose projection to 
$F$ (considered as a loop) is free homotopic to the solid loop in Figure~\ref{example2.fig}.
Let $K_1, K_2$ be Legendrian knots that are the same as $K$ everywhere except 
of a chart (contactomorphic to the standard contact $\R^3$) where $K_1$ and
$K_2$ are different from $K$ as it is described in
Figure~\ref{example4.fig}, see~\ref{description}. (Here $2r$ is the value of the Euler class of the
contact bundle on the homology class realized by $\mu(T^2)$.)

\begin{figure}[htbp]
 \begin{center}
  \epsfxsize 10cm
  \hepsffile{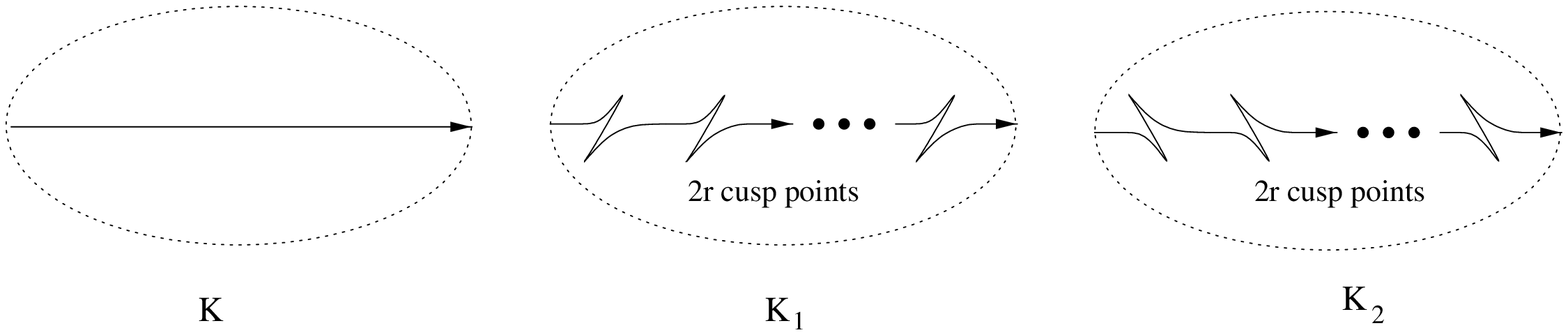}
 \end{center}
\caption{}\label{example4.fig}
\end{figure}

\end{emf}

\begin{thm}\label{example2} 
The knots $K_1$ and $K_2$ described above belong to the same
component $\mathcal L$ of the space of Legendrian curves and realize isotopic 
framed knots. There exists a $\Z$-valued order one invariant $I$ of 
Legendrian knots from $\mathcal L$ such that $I(K_1)\neq I(K_2)$.
\end{thm}

For the Proof of Theorem~\ref{example2} see Subsection~\ref{Proofexample2}.

\begin{rem} 
In Theorem~\ref{example2} we do not need the cooriented contact structure
$C$ to be tight or overtwisted. The only condition $C$ should satisfy is
the condition on the Euler class of the contact bundle. However the author
is not currently aware of the existence of tight contact satisfying the
conditions of Theorem~\ref{example2}.

Similar to~\ref{overtwisted} one can verify that the contact structure and
the knots $K_1$ and $K_2$ in the statement of Theorem~\ref{example2} can be
chosen so that the restrictions of the contact structure to the complements 
of $K_1$ and of $K_2$ are overtwisted.

Using the ideas of the Proof of Theorem~\ref{example2} one can 
easily construct many much more complicated examples of contact manifolds $(M,C)$ 
with a cooriented contact structure and connected components $\mathcal L$ of
the space of Legendrian curves in $(M, C)$ such that there exist Legendrian 
knots $K_1, K_2\in\mathcal L$ that realize isotopic framed knots but are distinguishable
by finite order invariants of Legendrian knots. For example instead of the
solid loop we could take any loop such that the number of double points that
separate the loop into two orientation reversing loops is odd, and the value
of the Euler class of the contact bundle on the homology class realized by
the mapping $T^2\rightarrow M$ corresponding to the loop is nonzero.
\end{rem}

\section{Transverse knots}
Let $M$ be a contact manifold with a
parallelized contact structure. Let $\mathcal T$ be a connected component of
the space of transverse curves in $M$, and let $\mathcal F$ be the connected
component of the space of framed curves that contains $\mathcal T$. (Such a
component exists because a transverse curve in a manifold with a
parallelized contact structure is naturally framed, and the path
in the space of transverse curves corresponds to a path in the space of
framed curves.) Let $V_n^{\mathcal T}$ (resp. $W_n^{\mathcal F}$) 
be the groups of $\mathcal A$-valued order $\leq n$ invariants of transverse
(resp. framed) knots from $\mathcal T$ (resp. from $\mathcal F$). Clearly every 
invariant $y\in W_n^{\mathcal F}$ restricted to
the category of transverse knots in $\mathcal T$ is an element $\phi(y)\in
V_n^{\mathcal T}$. This gives a homomorphism $\phi:W_n^{\mathcal
F}\rightarrow V_n^{\mathcal T}$.

\begin{thm}\label{isomorphismtransverse} Let $(M, C)$ be a contact manifold 
with a parallelized contact structure. Let $\mathcal T$ be a connected component of
the space of transverse curves in $(M,C)$, and let $\mathcal F$ be the
corresponding component of the space of framed curves.
Then the following two statements {\bf a:} and {\bf b:} are equivalent.
\begin{description}
\item[a] For any $x\in V_n^{\mathcal T}$ and any knots $K_1, K_2\in
\mathcal T$ that represent isotopic framed knots 
%(in $\mathcal F$)
$x(K_1)=x(K_2)$.
\item[b] 
For $x\in V_n^{\mathcal T}$ there exists $\psi (x)\in W_n^{\mathcal F}$ such 
that $\phi(\psi(x))=x$. Such $\psi (x)$ is unique and it defines a
canonical isomorphism $\psi:V_n^{\mathcal T}\rightarrow W_n^{\mathcal F}$.
\end{description}
\end{thm} 

The proof of Theorem~\ref{isomorphismtransverse} is analogous to the Proof 
of the Theorem~\ref{isomorphism}.

The following Theorem is an immediate consequence of
Theorems~\ref{thmTFtransverse} and~\ref{isomorphismtransverse}.

\begin{thm}\label{isomorphismobtainedtransverse}
Let $(M, C)$ be a contact manifold with a
parallelized contact structure, and let $\mathcal T$ be a connected component of
the space of transverse curves in $M$. Let $\mathcal F$ be the
connected component of the space of framed curves that contains $\mathcal T$.
Let $V_n^{\mathcal T}$ (resp. $W_n^{\mathcal F}$) 
be the groups of $\mathcal A$-valued order $\leq n$ invariants of transverse
(resp. framed) knots from $\mathcal T$ (resp. from $\mathcal F$). 
Then the groups $V_n^{\mathcal T}$ and $W_n^{\mathcal F}$ are canonically
isomorphic, provided that $\mathcal F$ satisfies condition \textrm{I}.
\end{thm}

\begin{rem}
Let $(M,C)$ be a contact manifold with a parallelized contact structure, and
let $\mathcal F$ be a connected component of the space of framed curves in
$M$ that satisfies condition \textrm{I}. Theorem~\ref{isomorphismobtainedtransverse} 
implies that for any $n\in\N$ the group of $\mathcal A$-valued order 
$\leq n$ invariants of transverse knots from a connected component 
of transverse curves contained in $\mathcal F$ does not depend on the
choice of a parallelized contact structure. And hence this group can not be
used to distinguish parallelized contact structures on $M$.
\end{rem}

\begin{rem}
Clearly if $\mathcal F$ consists of zero homologous knots, 
then it satisfies condition \textrm{I} and we get that for an arbitrary 
contact manifold $(M,C)$ (with a parallelized contact structure) 
and for an arbitrary component $\mathcal T\subset \mathcal F$
consisting of zero homologous transverse curves the groups $V_n^{\mathcal
T}$ and $W_n^{\mathcal F}$ are canonically   
isomorphic. 
Let $\mathcal T$ be a connected component of the space of
transverse curves in $\R^3$, and let $\mathcal F$ be the corresponding
component of the space of framed curves in $\R^3$.
Since $H_1(\R^3)=0$ and the standard contact structure in $\R^3$
is parallelized we get that the groups  $V_n^{\mathcal T}$
and $W_n^{\mathcal F}$ are canonically isomorphic. 
This generalizes the result of D.~Fuchs and S.Tabachnikov
~\cite{FuchsTabachnikov} that says that for the
standard contact $\R^3$ and for $\mathcal A=\C$ the quotient groups 
$V_n^{\mathcal T}/V_{n-1}^{\mathcal T}$ and $W_n^{\mathcal
F}/W_{n-1}^{\mathcal F}$ are canonically isomorphic. 
(There are two components of the space of transverse curves in $\R^3$. They
consist of respectively positive and negative transverse curves.
The connected components of the space of framed
curves in $\R^3$ are enumerated by the $\Z_2$-reduced value of the
self-linking number of knots from the component.)

\end{rem}

\begin{rem} For a big class of manifolds one can show that all the connected
components of the space of framed curves satisfy condition \textrm{I}. 
(See Theorems~\ref{fibration} and~\ref{atoroidal} and the remarks and
corollaries after them.) For an arbitrary parallelized contact structure on these 
manifolds one gets that the groups of $\mathcal A$-valued order $\leq n$ invariants 
of transverse and of framed knots (from the corresponding components) are canonically isomorphic.
\end{rem}

\section{Proofs}
\subsection{Useful Facts and Lemmas}
\begin{prop}\label{commute}
Let $p:X\rightarrow Y$ be a locally trivial $S^1$-fibration of an oriented
manifold $X$ over a (not necessarily orientable) manifold $Y$. Let
$f\in\pi_1(X)$ be the class of an oriented fiber of $p$, and let $\alpha$ be an
element of $\pi_1(X)$. Then:
\begin{description}
\item[a] $\alpha f=f\alpha\in\pi_1(X)$, provided that $p(\alpha)$ is an
orientation preserving loop in $Y$.
\item[b] $\alpha f=f^{-1}\alpha\in\pi_1(X)$, provided that $p(\alpha)$ is an
orientation reversing loop in $Y$.
\end{description}
\end{prop}

\begin{emf}{\em Proof of Proposition~\ref{commute}.\/}
If we move an oriented fiber along the loop $\alpha\in X$, then in the end it
comes to itself either with the same or with the opposite orientation. 
It is easy to see that it comes to itself with the opposite orientation if
and only if $p(\alpha)$ is an orientation reversing loop in $Y$.\qed
\end{emf}

\begin{prop}\label{Preissman}
Let $F\neq S^2, T^2\text{ (torus), } \R P^2, K\text{ (Klein bottle)}$
be a surface (not necessarily compact or orientable),
and let $G$ be a nontrivial commutative subgroup of $\pi_1(F)$.
Then $G$ is infinite cyclic.
%and there exists a
%unique maximal infinite cyclic $G'<\pi_1(F)$ containing $G$.
\end{prop}

\begin{emf}{\em Proof of Proposition~\ref{Preissman}.\/}
It is well known that any closed $F$, other than $S^2, T^2, \R P^2, K$,
admits a hyperbolic metric of a constant negative curvature.
(It is induced from the universal covering of $F$ by the hyperbolic plane
$H$.)
The Theorem by A.~Preissman (see~\cite{Docarmo} pp. 258-265)
says that if $M$ is a compact Riemannian manifold with a negative curvature,
then any nontrivial Abelian subgroup $G<\pi_1(M)$ is isomorphic to $\Z$.
Thus if $F\neq S^2, T^2, \R P^2, K$ is closed, then any nontrivial
commutative $G<\pi_1(F)$ is infinite cyclic.

%The proof of the Preissman's Theorem given in~\cite{Docarmo} is based on
%the fact, that if $\alpha,\beta\in\pi_1(M)$ are nontrivial commuting
%elements, then there exists a
%geodesic in $\bar M$ (the universal covering of $M$) that is mapped to
%itself
%under the action of these elements considered as deck transformations on
%$\bar M$. Moreover, these transformations restricted to the geodesic act as
%translations. This implies that if $F\neq S^2, T^2, \R P^2, K$
%is a closed surface, then there exists a unique maximal infinite cyclic
%$G'<\pi_1(F)$ containing $G$. This gives the proof of
%Proposition~\ref{Preissman} for
%closed $F$.

If $F$ is not closed, then the statement of the Proposition is also
true because in this case $F$ is homotopy equivalent to a bouquet of
circles.
\qed
\end{emf}

\begin{prop}\label{toughandtechnical}
Let $F\neq S^2, \R P^2, T^2, K$ (Klein bottle) be a surface
not necessarily closed or orientable. Let $M$ be an orientable manifold, and
let $p:M\rightarrow F$ be a locally trivial $S^1$-fibration. Let $f\in\pi_1(M)$
be the class of an oriented $S^1$-fiber of $p$, and let
$\alpha\in\pi_1(M)$ be an element with $p_*(\alpha)\neq 1\in\pi_1(F)$.
Let $\beta$ be an element from $Z(\alpha)<\pi_1(M)$ the centralizer of
$\alpha$.
Then:
\begin{description}
\item[a] If $p_*(\alpha)$ is an orientation preserving loop in $F$, then
there
exist $i,j\in\Z$ and $0\neq n\in\Z$ such that $\beta^n=\alpha^if^j$.
\item[b] If $p_*(\alpha)$ is an orientation reversing loop in $F$, then there
exist $i\in\Z$ and $0\neq n\in\Z$ such that $\beta^n=\alpha^i$.
\end{description}
\end{prop}

\begin{emf}\label{prooftoughandtechnical}
{\em Proof of Proposition~\ref{toughandtechnical}.\/}
Since $\beta\in Z(\alpha)$, we have that $\alpha$ and $\beta$ commute in
$\pi_1(M)$. Hence $p_*(\alpha)$ and $p_*(\beta)$ commute in $\pi_1(F)$.
Proposition~\ref{Preissman} and the fact that $p_*(\alpha)\neq 1\in\pi_1(F)$
imply that there exist $g\in\pi_1(M)$ with
$p_*(g)\neq 1 \in\pi_1(F)$, $i\in\Z$, and $0\neq n\in\Z$
such that $p_*(g)^n=p_*(\alpha)$ and $p_*(g)^i=p_*(\beta)$.

Hence,
\begin{equation}\label{equationtoughandtechnical}
\alpha=g^n f^k \text{ and } \beta=g^i f^l, \text{ for some } k,l\in\Z.
\end{equation}

Using Proposition~\ref{commute} we get that
\begin{equation}\label{equation2toughandtechnical}
\beta^{2n}=\alpha ^{2i} f^j,\text{ for some } j\in\Z.
\end{equation}
Since $n$ was initially chosen to be nonzero  we get
statement {\bf a:} of the proposition.

To prove statement {\bf b:} we have to
show that $j=0$ in~\eqref{equation2toughandtechnical}, provided that
$p_*(\alpha)$ is an orientation reversing loop
in $F$. Since $\beta\in Z(\alpha)$ we have that $\beta^{2n}$ commutes with
$\alpha$. Clearly $p_*(\alpha^{2i})$ is an orientation preserving loop in $F$,
and Proposition~\ref{commute} implies that $\alpha$ commutes with
$\beta^{2n}=\alpha^{2i}f^j$ if and only if $j=0$. (Note that $f$ has
infinite order in $\pi_1(M)$ for $M$ from the statement of the Proposition.)
\qed
\end{emf}

\begin{prop}\label{pi2nontrivial}
Let $M$ with $\pi_2(M)\neq 0$ be a connected three dimensional orientable manifold
that admits a structure of a locally-trivial $S^1$-fibration over a surface $F$ (not necessarily compact
or orientable). Then $M$ is either $S^1\times S^2$ or the quotient of
$S^1\times S^2$ by the equivalence relation $e^{2\pi i t}\times x=
e^{-2\pi i t}\times (-x)$. (In the last case $M$ fibers over $\R P^2$.)
\end{prop}

\begin{emf}{\em  Proof of Proposition~\ref{pi2nontrivial}.\/}
Let $p:M\rightarrow F$ be the locally trivial $S^1$-fibration. From the
exact homotopy sequence of the fibration we get that $F$ is either
$S^2$ or $\R P^2$.

Let $M$ be a manifold with $\pi_2(M)\neq 0$ that fibers over $S^2$. Choose
an orientation on $M$.
Since both $M$ and $S^2$ are oriented we get a canonical orientation on all
the $S^1$-fibers of $p$.
From the exact homotopy sequence of the fibration one gets that for such $M$
the class $f\in\pi_1(M)$ of the
fiber of $p$ has infinite order in $\pi_1(M)$. On the other hand if
$e\in\Z=H^2(S^2,\Z)$ is the Euler class of the oriented plane bundle associated 
to $p$, then $ef=0\in H_1(M)$. Since $\pi_1(M)$ is generated by $f$, we get that 
$f^e=1\in\pi_1(M)$, and hence $e=0$. But the only plane bundle over $S^2$ with
the zero Euler class is $\R^2\times S^2$ and it corresponds to 
$M=S^1\times S^2$.

Assume now that $M$ fibers over $\R P^2$. 
Let $p':M'\rightarrow S^2$ be the locally trivial $S^1$-fibration induced from 
$p:M\rightarrow \R P^2$ under the orientation double covering $S^2\rightarrow \R P^2$. 
Clearly $M'$ is orientable and $M'$ double covers $M$. Since $\pi_2(M)\neq
0$ we have that $\pi_2(M')\neq 0$, and hence $M'=S^1\times S^2$. Using this
one verifies that $M$ is 
the quotient of $S^1\times S^2$ by the equivalence relation $e^{2\pi i
t}\times x= e^{-2\pi i t}\times (-x)$. \qed
\end{emf}

\begin{emf}\label{Hansen} 
Let $X$ be a manifold, let $\Omega X$ be the space of free loops
in $X$, and let $\omega\in\Omega X$ be a loop. A loop $\alpha\in\pi_1
(\Omega X,\omega)$ is a mapping $\mu_{\alpha}:T^2=S^1\times S^1\rightarrow X$, 
with $\mu_{\alpha}(1, S^1)=\omega$ and $\mu_{\alpha}(S^1, 1)$ being the trace of the point 
$1\in S^1$ under the homotopy of $\omega$ described by $\alpha$. Let 
$t(\alpha)\in\pi_1(X, \omega(1))$ be the element corresponding to the trace
of the point $1\in S^1$ under the homotopy of $\omega$ described by $\alpha$.
Since $\pi_1(T^2)=\Z\oplus\Z$ is commutative, we get that 
$t:\pi_1(\Omega X,\omega)\rightarrow \pi_1(X, \omega(1))$ is a surjective
homomorphism of $\pi_1(\Omega X,\omega)$ onto the centralizer $Z(\omega)$
of $\omega\in\pi_1(X,\omega(1))$. 

If $t(\alpha)=t(\beta)\in\pi_1(X, \omega(1))$ for the loops 
$\alpha, \beta\in\pi_1(\Omega X, \omega)$, then the mappings $\mu_{\alpha}$
and $\mu_{\beta}$ of $T^2$ corresponding to these loops can be deformed to
be identical on the $1$-skeleton of $T^2$. Clearly the obstruction for 
$\mu_{\alpha}$ and $\mu_{\beta}$ to be homotopic as mappings of $T^2$ 
(with the mapping of the $1$-skeleton of $T^2$ fixed under homotopy) 
is an element of $\pi_2(X)$ obtained by gluing together the $2$-cells 
of the two tori. In particular we get the Proposition of V.~L.~Hansen~\cite{Hansen}) 
that $t:\pi_1(\Omega X,\omega)\rightarrow Z(\omega)<\pi_1(X,\omega(1))$ is an 
isomorphism, provided that $\pi_2(X)=0$.
\end{emf}

\begin{emf}\label{h-principleforcurves}
{\em $h$-principle for curves in $M$.}
For a three dimensional manifold $M$ we put $STM$ to be the manifold obtained by 
the fiberwise spherization of the tangent bundle of $M$, and we put
$\pr':STM\rightarrow M$ to be the corresponding locally trivial $S^2$-fibration.
The $h$-principle says that the space of curves in $M$ is weak homotopy
equivalent to $\Omega STM$ the space of free loops in $STM$. The weak
homotopy equivalence is given by mapping a curve $K$ to a loop $\vec
K\in\Omega STM$ that sends a point $t\in S^1$ to the point of $STM$
corresponding to the direction of the velocity vector of $K$ at $K(t)$.
\end{emf}

\begin{emf}\label{Proofexistcontact}{\em Proof of the
Proposition~\ref{existcontact}.\/} 
Since the contact structure is cooriented we get that the tangent bundle
$TM$ is isomorphic to a sum $C\oplus\epsilon$ of an oriented  contact bundle
$C$ and a trivial oriented line bundle $\epsilon$. 
But the tangent bundle of every
orientable $3$-manifold is trivializable and we get that the second Stiefel-Whitney 
class of the contact bundle is zero. But the second Stiefel-Whitney
of $C$ is the projection of the Euler class of $C$ under the natural mapping
$H^2(M, \Z)\rightarrow H^2(M, \Z_2)$ and we get the first statement of 
the Proposition.

Consider an oriented $2$-dimensional vector
bundle $\xi$ over $M$ with the Euler class $e(\xi)=e=2\alpha\in H^2(M,\Z)$.
Similar to the above we get that the second Stiefel-Whitney class $w_2(\xi)$ 
of $\xi$ is zero, since it is the projection of $e(\xi)=2\alpha\in H^2(M, \Z)$.
Since $\xi$ is an oriented vector bundle we have $w_1(\xi)=0$. 

Consider a sum $\xi\oplus\epsilon$ of $\xi$ with the trivial oriented
one-dimensional vector bundle $\epsilon$. Clearly the total Stiefel-Whitney class of the three
dimensional oriented vector bundle $\xi\oplus\epsilon$ is equal to $1$, and
the Euler class of $\xi\oplus\epsilon$ is equal to $0$. Using the
interpretation of the Stiefel-Whitney and the Euler classes of
$\xi\oplus\epsilon$
as obstructions for the trivialization of $\xi\oplus\epsilon$, we get that
$\xi\oplus\epsilon$ is trivializable. Since the tangent bundle of an oriented
three dimensional manifold is trivializable, we see that $\xi$ is isomorphic 
to an oriented sub-bundle of $TM$. Since $M$ is oriented this sub-bundle of
$TM$ is also cooriented.
Now the Theorem of Lutz~\cite{Lutz}, that says that every homotopy 
class of a distribution of tangent two-planes to $M$ contains a contact 
structure, implies the second statement of the Proposition.\qed
\end{emf}

\subsection{Proof of Theorem~\ref{isomorphism}}\label{proofisomorphism}
The fact that statement {\bf b} of Theorem~\ref{isomorphism} implies 
statement {\bf a} is clear.
Thus we have to  show that statement {\bf a} implies statement
\textrm{\bf b}.

Let $x\in V_n^{\mathcal L}$ be an order $\leq n$ invariant of Legendrian knots in
$\mathcal L$. In order to construct $\psi(x)\in W_n^{\mathcal F}$ we have to
specify the value of $\psi (x)$ on every framed knot from $\mathcal F$.

\begin{defin}[of $m(K_1, K_2)$ and of $K^0, K^{\pm 1}, K^{\pm
2}\dots$]\label{obstruction} 
Let $K_1$ and $K_2$ be two framed knots that 
coincide pointwise as embeddings of $S^1$. Then there is an integer obstruction 
$m(K_1, K_2)\in\Z$ for them to be isotopic as framed knots with the embeddings 
of $S^1$ fixed under the isotopy. This obstruction is calculated as follows.
Let $K'_1$ be the knot obtained by shifting $K_1$ along the framing and
reversing the orientation on the shifted copy. Together $K_1$ and $K'_1$
bound a thin strip. We put $m(K_1, K_2)$ to be the intersection number of 
the strip with the very small shift of $K_2$ along its framing. 

For a framed knot $K^0$ we denote by $K^i$, $i\in \Z$, the isotopy class of a
framed knot that coincides with $K^0$ as an embedding of $S^1$ and has
$m(K^0, K^i)=i$.

For two singular framed knots $K_{1s}$ and $K_{2s}$ with $n$ transverse 
double points
that coincide pointwise as immersions of $S^1$, we put $m(K_{1s},
K_{2s})\in\Z$ to be the value of $m$ on the nonsingular framed knots $K_1$ and $K_2$
that coincide pointwise as embeddings of $S^1$ and are obtained from
$K_{1s}$ and $K_{2s}$ by resolving each pair of the corresponding double points
of $K_{1s}$ and of $K_{2s}$ in the same way. (The value $m(K_{1s}, K_{2s})$
does not depend on the resolution as soon as the corresponding double points
of the two knots are resolved in exactly the same way.)
As before $m(K_{1s}, K_{2s})$ is the integer valued obstruction for $K_{1s}$
and $K_{2s}$ to be isotopic as singular framed knots with $n$ double points
with the immersion of $S^1$ corresponding to the two knots fixed under isotopy.

For a  singular framed knot $K^0_s$ with $n$ double  points we denote by $K^i_s$
the isotopy class of a singular framed knot with $n$ double points that
coincides with $K^0_s$ as an immersion of $S^1$ and has $m(K_s^0, K_s^i)=i$.
\end{defin}

\begin{prop}\label{homotopy} Let $K_1$ and $K_2$ be framed knots 
(resp. singular framed knots with $n$ double points) that coincide
pointwise as embeddings (resp. immersions) of $S^1$. 
Then $K_1$ and $K_2$ are homotopic as
framed knots (resp. singular framed knots with $n$ double points)
if and only if $m(K_1, K_2)$ is even.
\end{prop}

\begin{emf}\label{proofhomotopy}{\em Proof of Proposition~\ref{homotopy}.\/}
Clearly if $m(K_1, K_2)$ is even, then $K_1$ and $K_2$ are framed
homotopic. (We can change the obstruction by two by creating a small kink
and passing through a double point at it.)

Every oriented three-dimensional manifold $M$ is parallelizable. Hence all its 
Stiefel-Whitney classes vanish and $M$ admits a $\spin$-structure. A framed
curve $K$ in $M$ represents a loop in the principal $SO(3)$-bundle of $TM$. 
(The three-frame corresponding to a point of $K$ is the velocity
vector, the framing vector, and the unique third vector of unit length such
that the three-frame defines a positive orientation of $M$.)
One
observes that the values of the $\spin$-structure on the loops in the principal
$SO(3)$-bundle of $TM$ realized by $K_1$ and $K_2$ are different provided that 
$m(K_1, K_2)$ is odd. But these values do not change under homotopy of
framed curves. Hence if $m(K_1, K_2)$ is odd, then $K_1$ and $K_2$ are not
framed homotopic.
\qed
\end{emf} 
 
\begin{defin}[of $m(K)$]\label{strange} 
Using the self-linking invariant of framed knots one can easily
show that if $K_1$ and $K_2$ in~\ref{obstruction} are pointwise coinciding 
zero-homologous framed knots and  $m(K_1, K_2)\neq 0$,
then $K_1$ is not isotopic to $K_2$ in the category of framed knots. However in 
general this is not true, see~\ref{Proofexample1}.
This forces us to introduce the following definition.

If for an unframed knot $K$ there exist isotopic framed knots $K_1$ and $K_2$ 
that coincide with $K$ pointwise and have $m(K_1, K_2)\neq 0$, then we say that 
$K$ {\em admits finitely many framings\/}. 
For $K$ that admits finitely many
framings we put $m_K$ {\em the number of framings of\/} 
$K$ to be the minimal positive
integer $m$ such that there exist isotopic framed knots $K_1$ and $K_2$ that 
coincide with $K$ pointwise and have $m(K_1, K_2)=m$.
One can easily show that if $K$ admits finitely many framings, then there are
exactly $m_K$ isotopy classes of framed knots realizing the isotopy class of
the unframed knot $K$. Proposition~\ref{homotopy} implies that $m_K$ is even.

In a similar way we introduce the notion of the number of framings 
for unframed singular knots with $n$ double points.
\end{defin}

\begin{prop}\label{decrease} 
\begin{description}
\item[a] Let $K$ be an unframed knot obtained by forgetting a framing on a
knot from $\mathcal F$. Then there exists a Legendrian knot from $\mathcal
L$ realizing the isotopy class of $K$.
\item[b] If $K^0$ is an isotopy class of framed knots in $\mathcal F$ that 
is realizable by a Legendrian knot from $\mathcal L$, then the isotopy 
class of $K^{-2}$ is also realizable by a Legendrian knot from $\mathcal L$.
\item[c] Let $K_s$ be an unframed knot singular knot with $n$ double  points 
obtained by forgetting a framing on a singular knot  
from $\mathcal F$. Then there exists a Legendrian knot from $\mathcal  
L$ realizing the isotopy class of $K_s$.
\item[d] If $K_s^0$ is an isotopy class of framed knots in $\mathcal F$ that  
is realizable by a singular 
Legendrian knot from $\mathcal L$, then the isotopy class of $K_s^{-2}$ is
also realizable by a singular Legendrian knot from $\mathcal L$.
\end{description}
\end{prop}

\begin{emf}\label{proofdecrease}{\em Proof of Proposition~\ref{decrease}.\/}
{\em Proof of statement [a] of Proposition~\ref{decrease}.\/}
Let $CM$ be the fiberwise spherization of the two-dimensional
contact vector bundle, and let $\pr:CM\rightarrow M$ be the corresponding 
locally trivial $S^1$-fibration. We denote by $f\in \pi_1(M)$ the class of
an oriented $S^1$-fiber of $\pr$. 
For a Legendrian curve
$K_l:S^1\rightarrow M$ denote by $\vec K_l$ the loop in $CM$ obtained by
mapping a point $t\in S^1$ to the point of $CM$ corresponding to the
direction of the velocity vector of $K_l$ at $K_l(t)$. 

The $h$-principle proved for the Legendrian 
curve by M.~Gromov~\cite{Gromov} says that that $K_1$ and $K_2$ Legendrian
curves in $M$ belong to the same component of the space of Legendrian
curves in $M$ if and only if $\vec K_1$ and $\vec K_2$ are
free homotopic loops in $CM$. 

W.~L. Chow~\cite{Chow} and P.~K.~Rashevskii~\cite{Rashevskii} showed that
every knot $K$ is isotopic to a Legendrian knot $K_l$
(and this isotopy can be made $C^0$-small). 
Deforming $K$ we can assume (see~\ref{h-principleforcurves}) that: {\bf a:} $K$ and $K_l$ coincide in
the neighborhood of $1\in S^1\subset \C$, {\bf b:} $K$ and $K_l$ realize the
same element $[K]\in \pi_1(M, K_l(1))$, and {\bf c:} that liftings to $CM$ of 
Legendrian curves from $\mathcal L$ are free homotopic to a loop $\alpha$ in
$CM$ such that $\alpha(1)=\vec
K_l(1)$ and $\pr (\alpha)=[K] \in \pi_1(M, K_l(1))$. 

Proposition~\ref{commute} says that $f$ is in the
center of $\pi_1(CM, \vec K_l(1))$, since the contact structure is
cooriented and hence oriented. Then $\vec K_l=\alpha f^i\in \pi_1(CM, \vec K_l
(1))$, for some $i\in\Z$. 

Take a chart of $M$ (that is contactomorphic to the standard contact $\R^3$) 
containing a piece of the Legendrian knot.
From~\cite{FuchsTabachnikov} it is easy to see that insertions of 
two cusps shown in Figure~\ref{twocusp.fig}  
into the projection of a Legendrian knot in $\R^3$ onto the $(x,
z)$-plane induce multiplication by $f^{\pm 1}$ of the class in $\pi_1(CM,
\vec K_l(1))$ of a lifting of a knot $K_l$ to a loop in $CM$. (Here the sign depends on
the choice of an orientation of the fiber used to define $f$.) Performing
this operation sufficiently many times we obtain the Legendrian knot from
$\mathcal L$ realizing the unframed knot $K$.

\begin{figure}[htbp]
 \begin{center}
  \epsfxsize 10cm
  \hepsffile{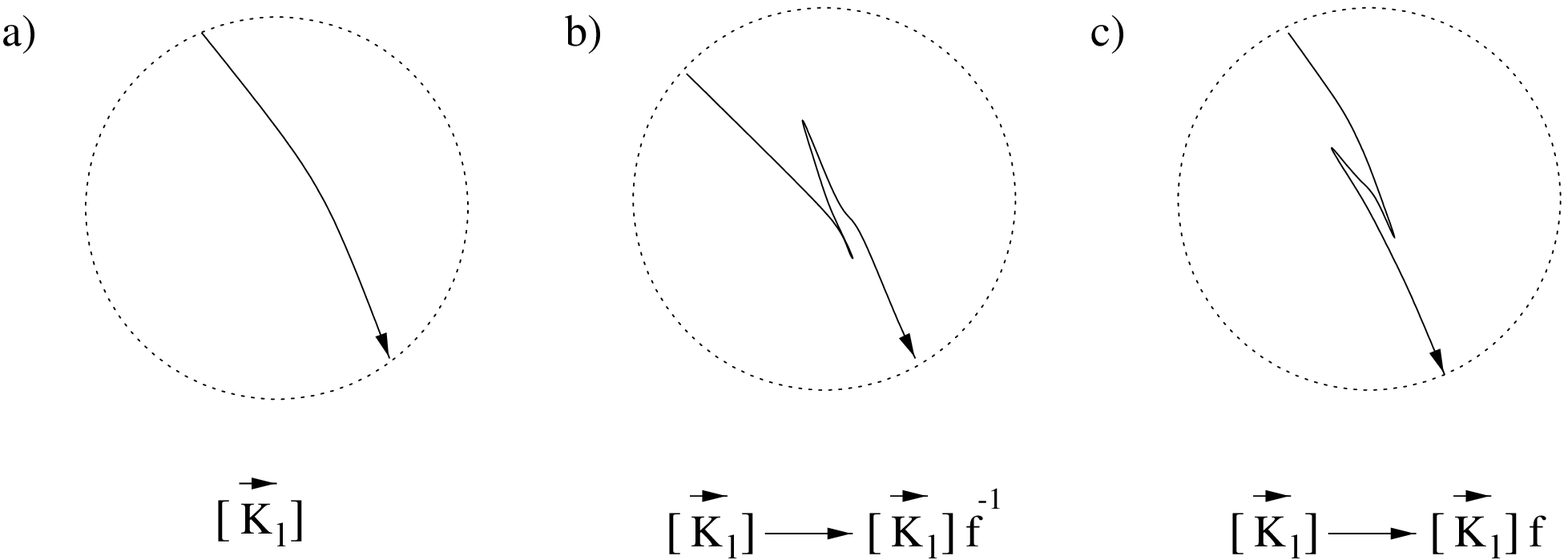}
 \end{center}
\caption{}\label{twocusp.fig}
\end{figure}

One easily modifies the arguments above to obtain 
the proof of statement [c] of Proposition~\ref{decrease}.

{\em Proof of statement [b] of Proposition~\ref{decrease}.\/}
Take a local chart of $M$ (where $M$ is contactomorphic to the standard contact
$\R^3$) containing a piece of the knot $K^0$ and perform the
homotopy in $\mathcal L$ shown in Figure~\ref{kink.fig}. (Observe that a
selftangency point of a projection of a Legendrian curve in $\R^3$ to the
$(x,z)$-plane corresponds to a double point of the Legendrian curve.)
A straightforward verification (see~\cite{FuchsTabachnikov}) shows that the Legendrian knot we 
obtain in the end of the homotopy realizes $K^{-2}$. 

\begin{figure}[htbp]
 \begin{center}
  \epsfxsize 12cm
  \hepsffile{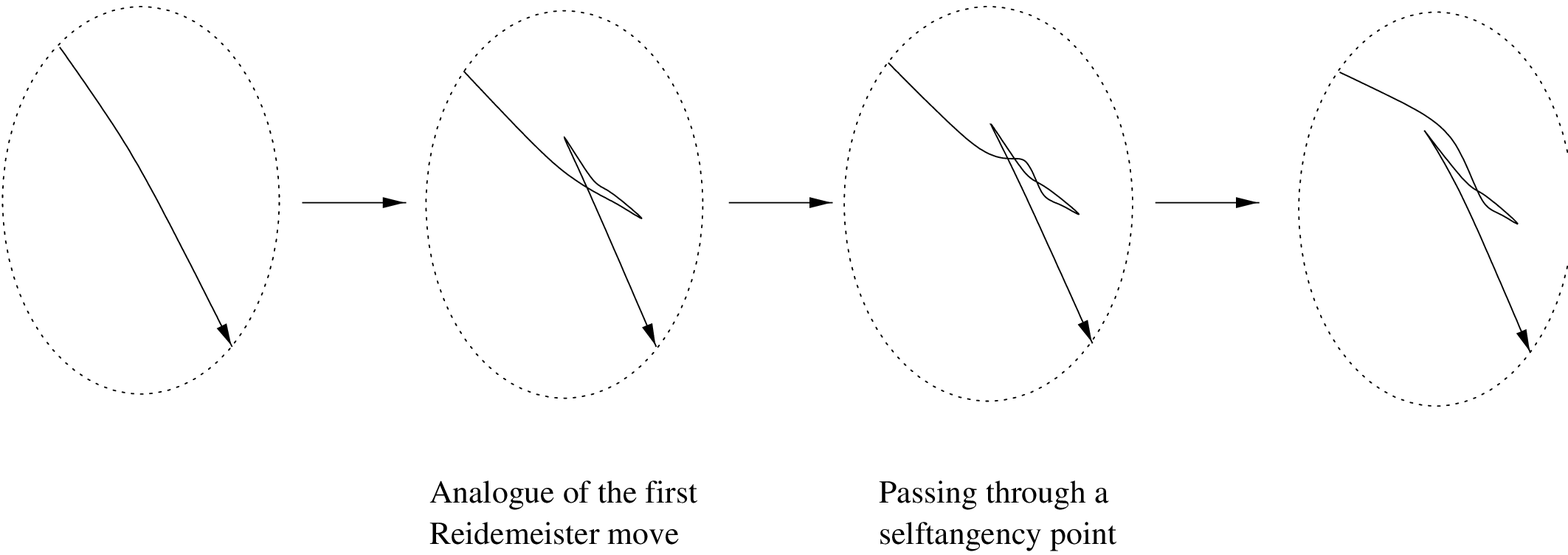}
 \end{center}
\caption{}\label{kink.fig}
\end{figure}

One easily modifies these arguments to obtain a Proof of Statement [d] of
Proposition~\ref{decrease}.
This  finishes the Proof of Proposition~\ref{decrease}.\qed 
\end{emf}

\begin{emf}{\em Definition of $\psi(x)$.\/} If an isotopy class of 
the knot $K\in \mathcal
F$ is realizable by a Legendrian knot $K_l\in\mathcal L$, then put
$\psi(x)(K)=x(K_l)$. The value $\psi(x)(K)$ is well defined because if
$K'_l\in\mathcal L$ is another knot realizing $K$, then
$x(K_l)=x(K'_l)$ by statement {\bf a:} of Theorem~\ref{isomorphism}.

Let $\mathcal C$ be the component of the space of unframed curves
that corresponds to forgetting framings on the curves from $\mathcal F$.
Propositions~\ref{decrease} and~\ref{homotopy} imply that 
if an unframed knot $K_u\in \mathcal C$
admits finitely many framings, then all the isotopy classes of framed knots
from $\mathcal F$ realizing the unframed knot $K_u$ are realizable by
Legendrian knots from $\mathcal L$. Thus we have defined the value of
$\psi(x)$ on all the framed knots from $\mathcal F$ that realize unframed 
knots admitting finitely many framings. 

If $K_u\in \mathcal C$ admits infinitely many framings, then either {\bf a)} all the
framed knots from $\mathcal F$ realizing the isotopy class of $K_u$ are
realizable by Legendrian knots from $\mathcal L$ or {\bf b)} there exists a
knot $K^0\in \mathcal F$  realizing the isotopy class of $K_u$ 
such that $K^0$ is realizable by a Legendrian knot
from $\mathcal L$ and $K^{+2}$ is not realizable by a Legendrian knot from 
$\mathcal L$ (in this case $K^{+4}, K^{+6}$ etc. also are not realizable by
Legendrian knots from $\mathcal L$, see~\ref{decrease}). In the case {\bf a)} the value of $\psi(x)$ is already defined for 
all the framed knots from $\mathcal F$ realizing $K_u$. In the case {\bf b)} put 
\begin{equation}\label{eqextension}
\psi(x)(K^{+2})=\sum_{i=1}^{n+1}\bigl((-1)^{i+1}\frac{(n+1)!}{i!(n+1-i)!}
\psi(x)(K^{+2-2i})\bigr).
\end{equation}
(Proposition~\ref{decrease} implies that the sum on the right hand side 
is well-defined.)
Similarly put 
\[
\begin{array}{l}

\psi(x)(K^{+4})=\sum_{i=1}^{n+1}\bigl((-1)^{i+1}\frac{(n+1)!}{i!(n+1-i)!}
\psi(x)(K^{+4-2i})\bigr), \\
\psi(x)(K^{+6})=\sum_{i=1}^{n+1}\bigl((-1)^{i+1}\frac{(n+1)!}{i!(n+1-i)!}
\psi(x)(K^{+6-2i})\bigr)\text{ etc.}\\
\end{array}
\]

Now we have defined $\psi(x)$ on all the framed knots (from $\mathcal F$) realizing $K_u$.
Doing this for all $K_u$ for which case {\bf b)} holds we define the value of $\psi(x)$
on all the knots from $\mathcal F$.

\end{emf}

{\bf Below we show that $\psi(x)$ is an order $\leq n$ invariant of framed knots from
$\mathcal F$.}
We start by proving the following Proposition.

\begin{prop}\label{mainidentity}
Let $K^0$ be a framed knot from $\mathcal F$, then $\psi(x)$ defined as
above satisfies identity~\eqref{eqextension}.
\end{prop}

\begin{emf}{\em Proof of Proposition~\ref{mainidentity}.\/}
If $K^{+2}$ is not realizable by a Legendrian knot from $\mathcal L$, then
the statement of the proposition follows from the formula we used to
define $\psi(x)(K^{+2})$. 

If $K^{+2}$ is realizable by a Legendrian knot $K_l$, then consider a
singular Legendrian knot $K_{ls}$ with $(n+1)$ double points that are
vertices of $(n+1)$ small kinks such that we get $K_l$
if we resolve all of them positively staying in the class of the Legendrian
knots. Let $\Sigma$ be the set of the $2^{n+1}$ ways in which we can
resolve $K_{ls}$. 
For $\sigma\in \Sigma$ put $\sign(\sigma)$ to be the sign of the resolution, and
put $K_{ls}^{\sigma}$ to be the nonsingular Legendrian knot obtained via the
resolution $\sigma$. Since $x$ is an order $\leq n$ invariant of Legendrian 
knots we get that
\begin{equation}
0=\sum_{\sigma\in \Sigma}\bigl(\sign(\sigma)x(K_{ls}^{\sigma})\bigr)=
\psi(x)(K_l)+\sum_{i=1}^{n+1}(-1)^i\frac{(n+1)!}{i!(n+1-i)!}\psi(x)(K_l^{-2i}).
\end{equation}
This finishes the proof of the Proposition. \qed
\end{emf}

\begin{emf}
Let $K_s\in\mathcal F$ be a singular framed knot with $(n+1)$ double points. Let $\Sigma$
be the set of $2^{n+1}$ ways of resolving the double points. For $\sigma\in
\Sigma$ put $\sign(\sigma)$ to be the sign of the resolution, and put $K^{\sigma}_s$
to be the isotopy class of the knot obtained via the resolution $\sigma$.

In order to prove that $\psi(x)$ is an order $\leq n$ invariant
of framed knots from $\mathcal F$, we have to show that 
\begin{equation}\label{identitytoshow}
0=\sum_{\sigma\in \Sigma}\bigl(\sign(\sigma)\psi(x)(K_{s}^{\sigma})\bigr),
\end{equation}
for every $K_s\in \mathcal F$.

First we observe that the fact whether the identity~\eqref{identitytoshow}
holds or not depends only on the isotopy class of the singular knot $K_s$ 
with $(n+1)$ double points. If the isotopy class of $K_s$ is realizable by 
a singular Legendrian knot with $(n+1)$ double points, then
identity~\eqref{identitytoshow} holds for $K_s$, since $x$ is an order $\leq n$
invariant of Legendrian knots and the value of $\psi(x)$ on a framed knot
$K$ realizable by a Legendrian knot $K_l$ was put to be $x(K_l)$.

Proposition~\ref{decrease} says that the singular unframed knot 
$K_{us}$ obtained by forgetting the framing on $K_s$ is realizable by a 
singular Legendrian knot from $\mathcal L$. 

If $K_{us}$ admits finitely many framings, 
then all the singular framed knots from $\mathcal F$ realizing the isotopy
class of $K_{us}$ are realizable by singular 
Legendrian knots from $\mathcal L$ and we get that
identity~\eqref{identitytoshow} holds for $K_s$. If $K_{us}$ admits infinitely 
many framings and all the singular framed knots from $\mathcal F$ realizing $K_{us}$ 
are realizable by singular Legendrian knots from $\mathcal L$, 
then~\eqref{identitytoshow} automatically holds for $K_s$. 

If $K_{us}$ admits infinitely many framings 
but not all the isotopy classes of singular framed knots from $\mathcal F$ 
realizing $K_{us}$ are 
realizable by singular Legendrian knots from $\mathcal L$. Then put $K_{us}^0$ to be the framed
knot realizing $K_{us}$ that is realizable by a singular Legendrian knot
from $\mathcal L$
and such that 
$K_{us}^{2i}$, $i>0$, are not realizable by singular Legendrian knots from
$\mathcal L$. 
Proposition~\ref{decrease} says that $K_{us}^{-2i}$, $i>0$, are realizable by
singular Legendrian knots from $\mathcal L$ and hence 
identity~\eqref{identitytoshow} holds for $K_{us}^{-2i}$, $i\geq 0$. 

Using Proposition~\ref{mainidentity} we 
show that~\eqref{identitytoshow} holds for $K_{us}^{+2}$. 
%\[
\begin{multline}
\sum_{\sigma\in
\Sigma}\sign(\sigma)\psi(x)({K_{us}^{+2}}^{\sigma})\\=
\sum_{\sigma\in
\Sigma}\Bigl(\sign(\sigma)
\sum_{i=1}^{n+1}(-1)^{i+1}\frac{(n+1)!}{i!(n+1-i)!}
\psi(x)({K_{us}^{(+2-2i)}}^{\sigma})\Bigr)\\ 
=
\sum_{i=1}^{n+1}\Bigl((-1)^{i+1}\frac{(n+1)!}{i!(n+1-i)!}\bigl(\sum_{\sigma\in
\Sigma}\sign(\sigma)\psi(x)({K_{us}^{(+2-2i)}}^{\sigma})\bigr)\Bigr)\\ 
=
\sum_{i=1}^{n+1}(-1)^{i+1}\frac{(n+1)!}{i!(n+1-i)!}\Bigl(0\Bigr)=0
\end{multline}
%\]

Similarly we show that ~\eqref{identitytoshow} holds for $K_{us}^{+4},
K_{us}^{+6}, \text{ etc}\dots$. 
\end{emf}

\begin{emf} Clearly $\psi$ is an injective homomorphism. 
Since an invariant $y\in W_n^{\mathcal F}$ should satisfy
identity~\eqref{eqextension} we get that $\phi$ is also injective. 
Finally we observe that $\phi\circ\psi=\id$ and this implies that 
$\psi:V_n^{\mathcal L}\rightarrow W_n^{\mathcal F}$ is an isomorphism.

To see that $\psi$ establishing the isomorphism is unique consider a 
singular framed knot with $(n+1)$ double points at the vertexes of $(n+1)$
small singular kinks. Since we want $\psi(x)$ to be an order $\leq n$ invariant, 
the sum of the values of it (with appropriate signs) on the knots obtained 
by the $2^{n+1}$ possible resolutions of the double points should be zero. This
forces us to use the formula~\eqref{eqextension} to define the value of 
$\psi(x)$ on the framed knots that are not realizable by Legendrian knots.

This finishes the Proof of Theorem~\ref{isomorphism}. \qed
\end{emf}

\subsection{Proof of Theorem~\ref{fibration}}\label{Prooffibration}
Assume that $\mathcal F$ does not satisfy condition \textrm{I}, then there
exists a framed knot $K_f\in\mathcal F$ such that $K_f$ is isotopic to
$K_f^{+2n}$ for some $0\neq n\in \Z$. (See~\ref{obstruction} for the
definition of $K_f^{+2n}$.) Let $\mathcal C$ be the connected component of
the space of unframed curves obtained by forgetting the framings on curves
from $\mathcal F$.

To a framed knot $K_f$ one can associate a link $L_f\in\mathcal C^2=\mathcal
C\times\mathcal C$ of two unframed knots with the first knot $K_1$ being $K_f$ 
with forgotten framing and the second knot $K_2$ being
obtained by a small shift of $K_1$ along the framing of $K_f$. 

Let $\alpha$ be a generic closed path in $\mathcal C^2$ starting with $L_f$.
Put $J_{\alpha}$ to be set of instances when the link becomes singular under the
deformation $\alpha$. At these instances $\alpha$ crosses the discriminant
in $\mathcal C^2$. (The discriminant is a
subspace of $\mathcal C^2$ formed by singular links.) Put $\sigma _j$, $j\in
J_{\alpha}$, to be the signs of these crossings.
Put $J'_{\alpha}\subset J_{\alpha}$ to be those instances when
the arising double point of the singular link is an intersection of the two
different components of the link. For a generic closed path $\alpha$
starting with $L_f$ put $\Delta_I(\alpha)=\sum_{j'\in J'_{\alpha}}\sigma _{j'}$.

Consider the instances when the first knot of the link becomes singular
under $\alpha$. The double point of $K_1$ corresponding to such an instance  
separates $K_1$ into two loops. Let $J''_{\alpha}\subset J_{\alpha}$ be the set of 
instances when both loops corresponding to a
double point of $K_1$ project to orientation reversing loops in $F$.
(If $F$ is orientable or if $\mathcal C$ consists of loops whose projections
to $F$ are orientation reversing loops, then $J''_{\alpha}=\emptyset$.) 
Similarly to the above put 
$\Delta_{I^o}(\alpha)=\sum_{j''\in J''}\sigma_{j^{''}}$.

For $\beta\in\pi_1(\mathcal C^2, L_f)$ put $t_1(\beta)\in \pi_1(M, K_1(1))$ to be
the trace in $M$ of the point of the knot $K_1$ corresponding to the point $1$
(on the unit complex circle that parameterizes $K_1$) under the deformation
described by $\beta$. It is easy to see that $t_1:\pi_1(\mathcal C^2,
L_f)\rightarrow \pi_1(M, K_1(1))$ is a homomorphism. 
Let $t_2:\pi_1(\mathcal C, K_2)\rightarrow \pi_1(M, K_2(1))$ be the
analogous homomorphism constructed using the trace in $M$ of the point on
$K_2$ corresponding to the point $1$ under the deformation of the link.

Clearly the groups $\pi_1(M, K_1(1))$
and $\pi_1(M, K_2(1))$ are naturally isomorphic. (For this reason in the
proof below we often identify these groups.)

Clearly if $K_f$ is isotopic to $K_f^{+2n}$, then there exists a generic closed path
$\alpha$ in $\mathcal C^2$ starting at $L_f$ such that $2n=\Delta_I(\alpha)$
and the restriction of 
$\alpha$ to each of the two components of the link is an isotopy.
Moreover this $\alpha$ should satisfy the condition that
$t_1(\alpha)=t_2(\alpha)\in \pi_1(M, K_1(1))=\pi_1(M, K_2(1))$.

\begin{emf}\label{idea}
To prove the Theorem we show that $\Delta_I(\alpha)=0$ 
for every generic closed $\alpha$ starting
with $L_f$  such that 
\begin{description}
\item[a]
$t_1(\alpha)=t_2(\alpha)\in\pi_1(M, K_1(1))$ and 
\item[b] 
$\alpha$ is realizable by
a deformation of $L_f$ the restriction of which to each of the two components of
the link is an isotopy.
\end{description}
\end{emf}

\begin{emf}\label{mainidea}
The codimension two stratum of the discriminant consists of links with two
distinct transverse double points. It is easy to see that
$\Delta_I(\beta)=0$ and $\Delta_{I^o}(\beta)=0$,
for any small loop $\beta$ going around the
codimension two stratum. This implies that $\Delta_I(\alpha)$ and 
$\Delta_{I^o}(\alpha)$
depend only on the element of $\pi_1(\mathcal C^2, L_f)$ realized 
by a generic loop $\alpha$.

Hence $\Delta_I$ and $\Delta_{I^o}$ are homomorphisms from $\pi_1(\mathcal C^2, L_f)$ 
to $\Z$. 
Clearly $\Delta_I(\alpha^p)=p\Delta_I(\alpha)$ 
and $\Delta_{I^o}(\alpha^p)=p\Delta_{I^o}(\alpha)$,
for every $\alpha\in\pi_1(\mathcal
C^2, L_f)$ and $p\in\Z$. 
Moreover if the class $\alpha$ satisfies conditions {\bf a:} and
{\bf b:} of~\ref{idea}, then $\alpha^p$ also satisfies these conditions.
Hence to prove the Theorem it suffices to show that for every
$\alpha\in\pi_1(\mathcal C^2, L_f)$ such that $t_1(\alpha)=t_2(\alpha)$
there exists $0\neq p\in\Z$ such that at least of the following holds:
\begin{description}
\item[1] $\Delta_I(\alpha^p)=0$, 
\item[2] $\Delta_{I^o}(\alpha ^p)\neq 0$ 
and hence $\alpha$ is not realizable by a deformation of $L_f$ the
restriction of which to every component of the link is an isotopy.
(The upper index $o$ in $I^o$ stands for obstruction.)
\end{description}
\end{emf}

\begin{emf}\label{obstructiondecrease} The $h$-principle (see~\ref{h-principleforcurves}) says that 
$\mathcal C$ is weak homotopy equivalent to $\Omega STM$, the space of free
loops in $STM$ the spherical tangent bundle of $M$. In particular
$\pi_1(\mathcal C, K_1)=\pi_1(\Omega STM, \vec K_1)$.
In Subsubsection~\ref{Hansen} we introduced a 
surjective homomorphism $t$ from $\pi_1(\Omega STM, \vec K_1)$ onto 
$Z(\vec K_1)<\pi_1(STM, \vec K_1(1))$. 

Let $\alpha, \beta\in\pi_1(\mathcal C, K_1)$ be loops such that $t(\alpha)=t(\beta)$. As it was explained 
in~\ref{Hansen} the obstruction for $\alpha$ and $\beta$ to be homotopic 
is an element of $\pi_2(STM)$. Since every closed orientable manifold is parallelizable 
we get that $STM=S^2\times M$. Proposition~\ref{pi2nontrivial} says that
$\pi_2(M)=0$ for $M$ from the statement of the Theorem~\ref{fibration}.
Hence $\pi_2(STM)=\pi_2(S^2)=\Z$. 

Let $k\in\Z$ be the obstruction for $\alpha$ and $\beta$ to be homotopic.
Consider the loop $\alpha'$ that looks the same as $\alpha$ except for a
small period of time when we perform the deformation shown in
Figure~\ref{obstruction.fig}. Clearly
$t(\alpha')=t(\alpha)=t(\beta)\in\pi_1(STM)$, and a straightforward
verification show that the obstruction for $\alpha'$ and $\beta$ to be
homotopic is $k-1$ (for an appropriate choice of the generator of
$\Z=\pi_2(STM)$). 
Performing this operation sufficiently many times we can change $\alpha$ 
to be homotopic to $\beta$. 

The kink participating in the deformation in
Figure~\ref{obstruction.fig} can be made very small, so that the deformation
happens far away from the second component of the link. Thus
$\Delta_I(\alpha)=\Delta_I(\alpha')$. Since one of the two
loops corresponding to a singular knot in Figure~\ref{obstruction.fig} is
contractible, it can not project to an orientation reversing loop in $F$.
Hence $\Delta_{I^o}(\alpha)=\Delta_{I^o}(\alpha')$. 

Clearly $\pi_1(STM)=\pi_1(M)$ and 
\begin{equation}
Z(\vec K_1)=Z(\vec K_2)=Z(K_1)=Z(K_2)<
\pi_1(M, K_1(1))=\pi_1(M, K_2(1)). 
\end{equation}
(Here the equalities denote the canonical
isomorphisms and $Z(K_1)$ (resp. $Z(K_2)$) is the centralizer of $K_1$
(resp. $K_2$) in $\pi_1(M, K_1(1))$.)
Hence the above observations imply that for $\gamma\in \pi_1(\mathcal C^2, L_f)$ 
such that $t_1(\gamma)=t_2(\gamma)\in Z(K_1)=Z(K_2)<\pi_1(M, K_1(1))=\pi_1(M, K_2(1))$
the quantities $\Delta_I(\gamma)$ and $\Delta_{I^o}(\gamma)$ depend only on
the element of $Z(K_1)<\pi_1(M, K_1(1))$ realized by $t_1(\gamma)$.

\begin{figure}[htbp]
 \begin{center}
  \epsfxsize 10cm
  \hepsffile{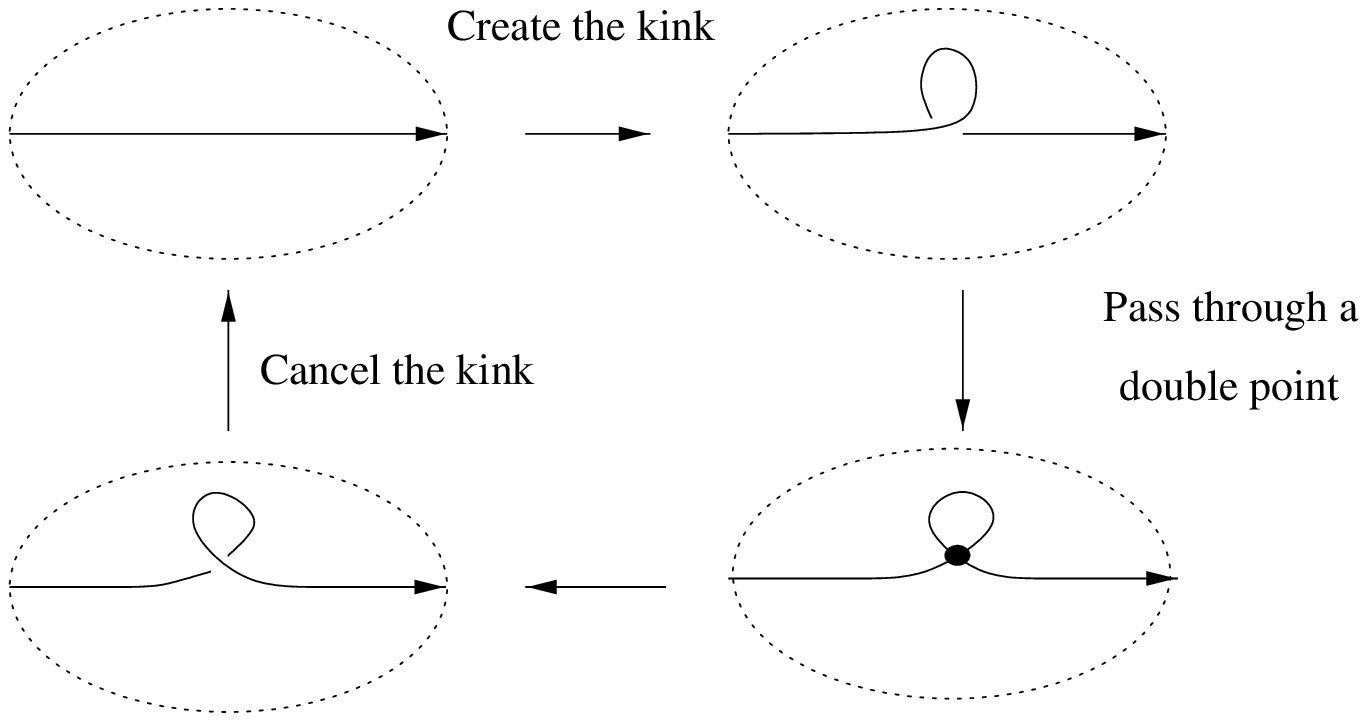}
 \end{center}
\caption{}\label{obstruction.fig}
\end{figure}
\end{emf}

\begin{emf}\label{mainidea2}
Combining this with~\ref{mainidea} we get that to prove the Theorem it
suffices to show that for every
$\beta\in Z(K_1)$ there exist $0\neq n\in\Z$ and $\gamma\in\pi_1(\mathcal
C^2, L_f)$ such that $t_1(\gamma)=t_2(\gamma)=\beta^n\in \pi_1(M, K(1))$ and
either $\Delta_I(\gamma)=0$ or $\Delta_{I^o}(\gamma)\neq 0$. 
\end{emf}

Let $p:M\rightarrow F$ be a fibration from the statement of the
proposition, and let $f\in\pi_1(M, K_1(1))$ be the class of an oriented 
$S^1$-fiber of $p$. 

{\bf Proof of the Theorem for $F\neq S^2, T^2, \R
P^2, K$.\/}
There are two cases either $p_*(K)=1\in\pi_1(F)$ or not. 
For $F\neq S^2, T^2, \R  P^2, K$ we prove the Theorem in these two cases
respectively in~\ref{notcontractible} and~\ref{contractible}.

\begin{emf}\label{notcontractible}
{\em Consider the case of $p_*(K)\neq 1\in\pi_1(F)$ for $F\neq S^2, T^2, \R   
P^2, K$.\/} 
Proposition~\ref{toughandtechnical} says that if $p(K_1)$ is an
orientation reversing loop in $F$, then for every $\beta\in
Z(K_1)$ there exist $0\neq n\in\Z$ and $i\in\Z$ such that $\beta^n=K_1^i$.
The loop $\gamma\in\pi_1(\mathcal C^2, L_f)$ that has
$t_1(\gamma)=t_2(\gamma)=\beta^n$ can be realized as 
the $i$-th power of 
$\gamma_2$, for the {\em loop \/}$\gamma_2\in\pi_1(\mathcal C^2, L_f)$ 
that is the sliding of both knots
$K_1$ and $K_2$ along themselves according to their orientations. 
(This deformation is induced by the rotation of the circle that parameterizes 
$K_1$ and $K_2$.)
Clearly no singular links occur under the deformation described by
$\gamma_2$. Hence $\Delta_{I}(\gamma_2^i)=i\Delta_{I}(\gamma_2)=0$ and we
have proved the statement of the Theorem in the case of $p(K_1)$ being an 
orientation reversing loop in $F\neq S^2, \R P^2, T^2, K$.

Proposition~\ref{toughandtechnical} says that if $p(K_1)$ is an
orientation preserving loop in $F$, then for every $\beta\in
Z(K_1)$ there exist $0\neq n\in\Z$ and $i,j\in\Z$ 
such that $\beta^n=K_1^if^j$. The loop $\gamma\in\pi_1(\mathcal C^2, L_f)$
that has $t_1(\gamma)=t_2(\gamma)=\beta^n$ can be realized as a product of the 
$j$-th power of the loop $\gamma_1$ and the $i$-th power of the loop
$\gamma_2$, for $\gamma_1$ described below (and for $\gamma_2$ described
above).

{\em Loop $\gamma_1$.\/} Since $p(K_1)$ is an orientation preserving loop 
and $M$ is orientable, we get that the $S^1$-fibration over $S^1$ 
(parameterizing the knots) induced by $p\circ K_1:S^1\rightarrow F$ is
trivializable. 
Hence we can coherently orient the fibers of this fibration. The orientation
of the $S^1$-fiber over $t\in S^1$ induces the orientation of the $S^1$-fiber
of $p$ that contains $K_1(t)$. 
Similarly for each $t\in S^1$ we have an induced orientation of the 
$S^1$-fiber of $p$ that contains $K_2(t)$. (In these two constructions we
use the same orientation of the $S^1$-fibers of $S^1\times S^1\rightarrow S^1$ 
so that for every $t\in S^1$ the orientations of the fibers containing $K_1(t)$ 
and $K_2(t)$ are coherent.) 
The loop $\gamma_1$ is the deformation of $L_f$ 
under which every point of $L_f$ slides once around the fiber of $p$  
(staying inside the fiber that contains the point) in the
direction specified by the orientation of the fiber.

As it was explained in~\ref{mainidea2} to prove the Theorem 
it suffices to show that 
either $\Delta_I(\gamma)=0$ or
$\Delta_{I^o}(\gamma)\neq 0$.
We show that 
\begin{equation}\label{2IO=I}
2\Delta_{I^o}(\gamma)=
\Delta_I(\gamma).
\end{equation} 
Clearly this gives the proof of the Theorem.

The input into $\Delta_I(\gamma_1)$ could come only from singular links
whose double points are located over the crossing points of $p(K_1)\cap p(K_2)$.
Similarly the input into $\Delta_{I^o}(\gamma_1)$ could come only from the 
double points of $K_1$ (arising under the deformation $\gamma_1$) 
and they could occur only over the selfintersection points of $p(K_1)$.

The important crossing points of $p(K_1)\cap p(K_2)$ are located in pairs in
small neighborhoods of the selfintersection points of $p(K_1)$, see
Figure~\ref{doublepoint.fig}. (Other
points of $p(K_1)\cap p(K_2)$ correspond to extra twists of the framing. 
The two orientations of the fiber over them are the same, and hence the
movement of branches of $K_1$ and of $K_2$ over such points happen in the same
direction and does not give rise to singular links.)

If the double point $x$ of $p(K_1)$ corresponding to a pair $(x', x'')$ 
of crossing points of $p(K_1)\cap p(K_2)$ separates $p(K_1)$ into 
two orientation preserving loops, then the sliding of the branches of 
$K_1$ and of $K_2$ located over $x'$ and $x''$ happen in the same direction
and no double points of $K_1\cap K_2$ occur under this process. Hence 
such $x$ corresponds to the zero input into $\Delta_I(\gamma_1)$. On the
other hand the sliding of the two branches of $K_1$ located over such $x$ also
happens in the same direction. Thus the input into $\Delta_{I^o}(\gamma_1)$
corresponding to such $x$ is also zero.

\begin{figure}[htbp]
 \begin{center}
  \epsfxsize 12cm
  \hepsffile{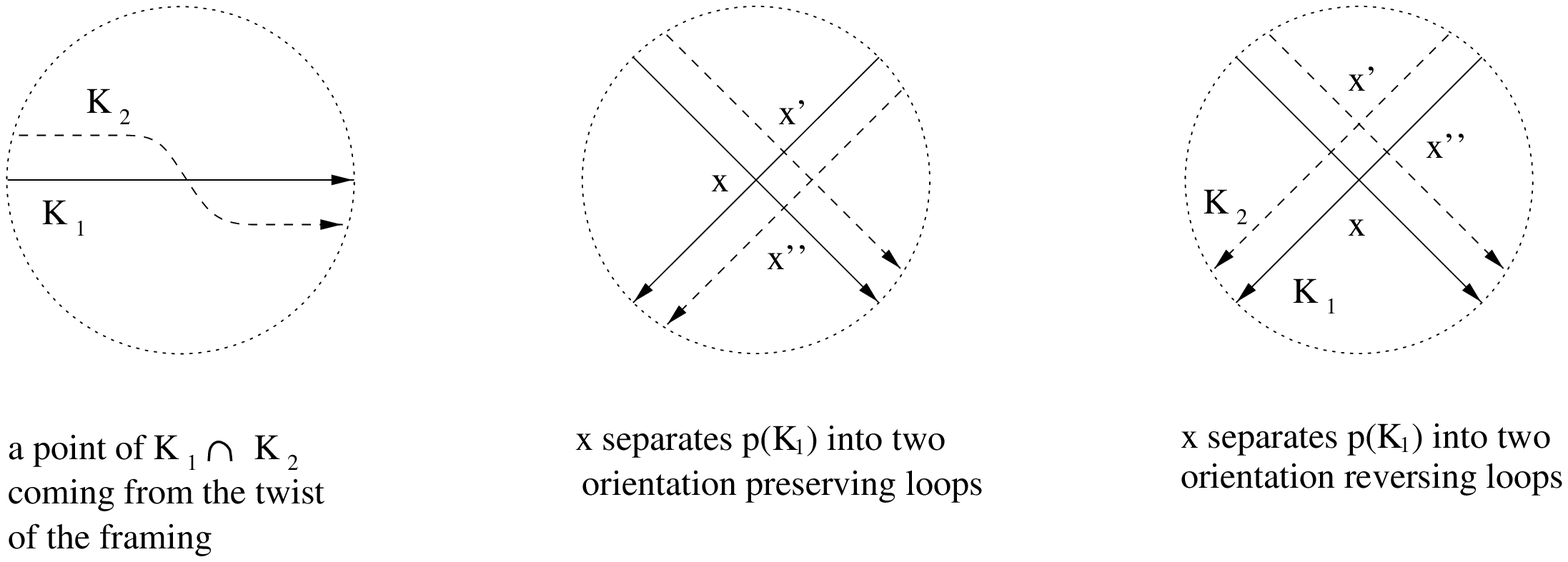}
 \end{center}
\caption{}\label{doublepoint.fig}
\end{figure}

If the double point $x$ of $p(K_1)$ corresponding to a pair $(x', x'')$ 
of crossing points of $p(K_1)\cap (K_2)$ separates $p(K_1)$ into 
two orientation reversing loops, then the sliding of the branches of $K_1$
and of $K_2$ located over $x'$ (and over $x''$) happens in the opposite directions 
and such $x$ corresponds to the nontrivial input into $\Delta_I(\gamma_1)$.
On the other hand the sliding of the two branches of $K_1$ located
over $x$ also happens in the opposite directions, and such $x$ corresponds 
to a nontrivial input $\Delta_{I^o}(\gamma_1)$.
A straightforward verification show that for such $x$ 
the input into $\Delta_I(\gamma_1)$
corresponding to $x$ is twice the input into $\Delta_{I^o}(\gamma_1)$
corresponding to $x$.

Thus we have $\Delta_{I}(\gamma_1)=2\Delta_{I^o}(\gamma_1)$. Clearly 
$\Delta_{I}(\gamma_2)=\Delta_{I^o}(\gamma_2)=0$. 
Substituting these expressions into 
$\Delta_{I^o}(\gamma)=j\Delta_{I^o}(\gamma_1)+ i\Delta_{I^o}(\gamma_2)$ and into 
$\Delta_{I}(\gamma)=j\Delta_{I}(\gamma_1)+i\Delta_{I}(\gamma_2)$ 
we obtain identity~\eqref{2IO=I}. This finishes the Proof of
Theorem~\ref{fibration} for $p(K_1)$ being not contractible in $F\neq S^2,
\R P^2, T^2, K$.
\end{emf}

\begin{emf}\label{contractible}
{\em Consider the case of $p_*(K_1)=1\in\pi_1(F)$ and $F\neq S^2, \R
P^2, T^2, K$.\/} In this case $K_1=K_2=f^n\in\pi_1(M)$. For all these
surfaces $f$ has infinite order in $\pi_1(M)$. Proposition~\ref{commute}
implies that if $n\neq 0$, then $Z(K_1)$ is $\pi_1^{pres}(M, K_1(1))$, the
subgroup of $\pi_1(M)$ consisting of the elements whose projection to $F$
are orientation preserving loops in $F$. 

If $n=0$ (and hence $K_1=1\in\pi_1(M)$), then $Z(K_1)=\pi_1(M)$. Clearly in
this case $\alpha^2\in \pi_1^{pres}(M, K_1(1))$, for every 
$\alpha\in\pi_1(M)=Z(K_1)$. 

From~\ref{mainidea2} we see that in both these cases to prove the Theorem it
suffices to show that for every $\beta\in\pi_1^{pres}(M, K_1(1))$ there
exists $\gamma\in\pi_1(\mathcal C^2,L_f)$ with
$t_1(\gamma)=t_2(\gamma)=\beta$ and $\Delta_I(\gamma)=0$.

Such $\gamma$ can be realized as follows.
Consider a homotopy $r$ that changes $K_1$ and $K_2$ so that $p(K_1)$
and $p(K_2)$ are located in two small disks that do not intersect. 
(We can assume that $K_1(1)$ and $K_2(1)$ do not move under this homotopy.) 
Next consider the deformation of the link described by the loop $\alpha$ in 
$\mathcal C^2$ such that $t_1(\alpha)=\beta$, $t_2(\alpha)=\beta$, and at each moment of 
the deformation $p(K_1)$ and $p(K_2)$ are in small disks that do not intersect. 
(The position of the disks in $F$ changes under $\alpha$.) 
Finally deform $K_1$ and $K_2$ back to the original
shape via $r^{-1}$. Clearly $\Delta_I(\alpha)=0$ and the inputs of $r$ and
$r^{-1}$ into $\Delta_I(\gamma)$ cancel out. Thus $\Delta_I(\gamma)=0$ and
we finished the proof of the Theorem in the case of $F\neq S^2, \R P^2, T^2, K$.
\end{emf}

\begin{emf} {\em Proof of Theorem~\ref{fibration} for $F=S^2$.\/}
In this case $\pi_1(M)$ is generated by the class $f$ of the fiber of
$p:M\rightarrow S^2$. Proposition~\ref{pi2nontrivial} says that $\pi_2(M)$
is zero (for $M$ from the statement of the Theorem). Hence $f$ has finite
order in $\pi_1(M)$, and a certain degree of every $\beta\in\pi_1(M, K_1(1))$
is $1\in\pi_1(M)$. The loop $1\in\pi_1(M)$ is realizable as $t_1(\gamma)=t_2(\gamma)$ for
a trivial deformation $\gamma$ (under which nothing happens and hence 
no singular link occur). This finishes the proof of the Theorem for $F=S^2$.
\end{emf}

\begin{emf} {\em Proof of Theorem~\ref{fibration} for $F=\R P^2$.\/}
Proposition~\ref{pi2nontrivial} says that $\pi_2(M)$   
is zero (for $M$ from the statement of the Theorem). Hence $f$ has finite
order $k$ in $\pi_1(M)$. Using Proposition~\ref{commute} and the fact that 
$\pi_1(\R P^2)=\Z_2$ we get that $\beta^{2k}=1\in\pi_1(M)$, for every
$\beta\in\pi_1(M)$.
The loop $1\in\pi_1(M)$ is realizable as $t_1(\gamma)=t_2(\gamma)$ for
a trivial deformation $\gamma$ (under which nothing happens and hence 
no singular link occur) and we have proved the Theorem for $F = \R P^2$. 
\end{emf}

\begin{emf}{\em Proof of Theorem~\ref{fibration} for $F=T^2$.\/}
Let $p:M\rightarrow T^2$ be the fibration. Since $M$ is oriented and the class
$f$ of the fiber of $p$ is in the center of $\pi_1(M)$ (see~\ref{commute})
we get that 

\begin{equation}\label{group}
\pi_1(M)=\bigl\{a, b\big| a b a^{-1} b^{-1}=f^k, af=fa, bf=fb\bigr\}, \text{
for some } k\in\Z.
\end{equation}
(Here $a$ and $b$ are elements projecting respectively to the meridian and 
longitude of $T^2$, and $k\in\Z=H^2(T^2, \Z)$ is the Euler class of $p:M\rightarrow T^2$.) 

From~\eqref{group} one gets that 
\begin{equation}\label{formula}
b^i a^j=f^{-kij} a^j b^i, \text{ for all } i,j\in\Z.
\end{equation}
Every element of $\pi_1(M)$ can be uniquely presented as 
$a^l b^m f^n$, for some $l,m,n\in\Z$. Let $K_1=a^l b^m f^n\in\pi_1(M)$. 
Using~\eqref{formula} and the fact that $f$ is in the center of $\pi_1(M)$ one 
verifies that $a^{l_1}b^{m_1}f^{n_1}$ commutes with $a^lb^mf^n$ if and only if
$kml_1=km_1l$. From this one concludes that if $k\neq 0$ and $p(K_1)\neq 1\in\pi_1(T^2)$, 
then a certain degree of every element of $Z(K)$ is expressible as
$t_1(\gamma)=t_2(\gamma)$ for $\gamma$ that is a product
of the powers of the loops $\gamma_1$ and $\gamma_2$ described
in~\ref{notcontractible}. Similarly to~\ref{notcontractible} we get the
proof of the Theorem for this case. 
If $p(K)=1\in\pi_1(T^2)$ and $k\neq 0$, then the proof of the Theorem is 
analogous to~\ref{contractible}.

If $k=0$, then $M=S^1\times S^1 \times S^1$ and $M$ has three structures of
an $S^1$-fibration over $T^2$. (They are obtained by projecting $M$ on the three
different pairs of circles.) Then every element of $\pi_1(M)$ is realizable
as $t_1(\gamma)=t_2(\gamma)$ for $\gamma$ being a product of loops that are the
powers of the versions of the loop $\gamma_1$ with respect to the three
$S^1$-fibration structures. We get that $\Delta_I(\gamma)=0$ and this
finishes the proof of Theorem~\ref{fibration} in the case of $F=T^2$.
\end{emf}

\begin{emf}\label{K}{\em Proof of Theorem~\ref{fibration} in the case of 
$F=K$.\/}
Using Proposition~\ref{commute}, the fact that the class $f\in\pi_1(M)$ 
of the $S^1$-fiber of $p:M\rightarrow K$ is of infinite order in $\pi_1(M)$, 
and the fact that a relation on the generators of
$\pi_1(M)$ projects (under $p_*$) to a relation on the generators of
$\pi_1(K)$ we obtain that 
\begin{equation}\label{pi1MK}
\pi_1(M)=\bigl \{ g,h,f\big| hf=f^{-1}h, gf=fg, hg=g^{-1}hf^{k}\}, 
\text{ for some } k\in\Z .
\end{equation}
(Here $g$ and $h$ project respectively to $c$ and $d$ in
Figure~\ref{klein1.fig} depicting the Klein bottle.)

\begin{figure}[htbp]
 \begin{center}
  \epsfxsize 3cm
  \hepsffile{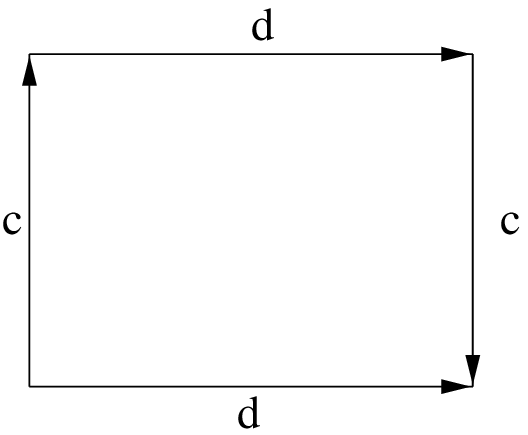}
 \end{center}
\caption{}\label{klein1.fig}
\end{figure}

Using~\eqref{pi1MK} one verifies that:
\begin{equation}\label{multtablepi1MK} 
h^ig^j=
\begin{cases}
g^{-j}h^if^{kij}& i \text{ odd},\\
g^jh^if^{kij}& i \text{ even}.
\end{cases}
\end{equation}

Every element of $\pi_1(M)$ can be uniquely presented as $g^lh^mf^n$,
for some $l,m,n\in\Z$. Taking $K_1=g^lh^mf^n\in\pi_1(M)$ and
$\beta=g^{l_1}h^{m_1}f^{n_1}\in\pi_1(M)$ we find all the values of
$l,m,n,l_1,m_1,n_1\in\Z$ such that $\beta\in Z(K_1)<\pi_1(M)$.
We obtain that $\beta \in Z(K_1)$ if and only if one of the following
holds:
\begin{description}
\item[1] $m$ and $m_1$ are even, and $kl_1m=klm_1$;
\item[2] $m$ is even, $m_1$ is odd, $l=0$, and $2n=-kl_1m$;
\item[3] $m$ is odd, $m_1$ is even, $l_1=0$, and $2n_1=-klm_1$;
\item[4] $m$ and $m_1$ are odd, $l_1=l$, and $2(n-n_1)=kl(m_1-m)$.
\end{description}

A straightforward verification shows that: 
\begin{description} 
\item[a]if $m$ is odd, then a certain
nonzero power of every $\beta\in Z(K)$ can be realized as 
$t_1(\gamma_2^q)=t_2(\gamma_2^q)$,
for some $q\in\Z$ (see~\ref{notcontractible} for the definition of
$\gamma_2$).
And similarly to~\ref{notcontractible} we get the proof of the Theorem
for this case.
\item[b] if $m$ is even and $k$ in~\eqref{pi1MK} is nonzero, then a certain
nonzero power of every $\beta\in Z(K)$ can be realized as 
$t_1(\gamma_1^p\gamma_2^q)=t_2(\gamma_1^p\gamma_2^q)$, for some $p, q\in\Z$ 
(see~\ref{notcontractible} for the definitions of $\gamma_1, \gamma_2$).
And similarly to~\ref{notcontractible} we get the proof  for this case.
\item[c] if $m$ is even, $k=0$, and $K_1=h^m$, then $Z(K)$ is generated by 
$g,h,f$.
\item[d] if $m$ is even, $k=0$, and $K_1=g^lh^mf^n$ with at least one of $l$ and $n$
being nonzero, then $Z(K)$ is generated by $g, h^2, f$. 
\end{description}

To finish the proof of Theorem~\ref{fibration} we consider the cases {\bf c}
and {\bf d}. 
In both cases a certain nonzero power of every
$\beta\in Z(K_1)$ belongs to the subgroup of $\pi_1(M)$ generated by 
$g, h^2, f$. 

Consider an $S^1$-fibration $p':M'\rightarrow T^2$ induced by the
orientation double cover $T^2\rightarrow K$. Clearly $M'$ is orientable and 
using~\eqref{pi1MK} we get that $\pi_1(M')=\{ g, h^2, f| g h^{2}= h^2 g, gf=fg,
h^2f=fh^2\}$. Using the fact that $k$ in~\eqref{group} is the Euler class of the
$S^1$-fibration over $T^2$, we see that $M'=S^1\times S^1 \times S^1=STT^2$. 
This implies that either $M$ is $STK$ (the spherical tangent bundle of the
Klein bottle) or the fiberwise projectivization of $M$ is $STK$.

We give the proof in the case of $M=STK$. The proof in the case when the
fiberwise projectivization of $M$ is $STK$ is obtained is the similar way.
  
Consider a loop $\mu$ in the space of autodiffeomorphisms of $STK$ that is
induced by the sliding of $K$ along the unit vector field parallel to the
curve $d$ on $K$ (see Figure~\ref{klein1.fig}). (Note that $K$ has to slide twice along itself under this
loop before every point of it comes to the original position.) This $\mu$
induces an isotopy $\gamma_3$ of $L_f$ such that
$t_1(\gamma_3)=t_2(\gamma_3)=h^2f^i$, for some $i\in\Z$. Since
$\gamma_3$ is an isotopy we have that
$\Delta_{I}(\gamma_3)=\Delta_{I^o}(\gamma_3)=0$. As usually 
$t_1(\gamma_1)=t_2(\gamma_1)=f$
and $\Delta_{I}(\gamma_1)=2\Delta_{I^o}(\gamma_1)$. Next we observe that
$STK$ has another structure of the $S^1$-fibration over $K$ for which $g$ is
the class of the $S^1$-fiber and the elements $h,f$ project to the generators of
$\pi_1(K)$. Let $\gamma_4\in\pi_1(\mathcal C^2, L_f)$ be the loop that is
the version of $\gamma_1$ with respect to this new structure of the
$S^1$-fibration on $STK$. Clearly $t_1(\gamma_4)=t_2(\gamma_4)=g$ and
$\Delta_{I}(\gamma_4)=2\Delta_{I^o}(\gamma_4)$. Similarly
to~\ref{notcontractible} this finishes the proof of Theorem~\ref{fibration} in 
the case when $M=STK$.

This finishes the Proof of Theorem~\ref{fibration} for all the surfaces.\qed
\end{emf}

\subsection{Proof of Proposition~\ref{interpretationconditionII}.\/}
\label{ProofinterpretationconditionII}
The $h$-principle for curves~\ref{h-principleforcurves} says
that the set $\mathcal C$ 
of the connected components of the space of curves in $M$ is
naturally identified with the set of the connected components of the space of
free loops in $STM$ the spherical tangent bundle of $M$, or which is the
same with the set of conjugacy classes of the elements of $\pi_1(STM)$. From
the long homotopy sequence of the fibration $\pr':STM\rightarrow M$ we see
that it is also naturally identified with the set of conjugacy classes of
the elements of $\pi_1(M)$. Choose a $\spin$-structure on $M$. It is easy to
see (cf.~\ref{homotopy} and~\ref{proofhomotopy}) that the set 
$\mathcal C_{\mathcal F}$ of the connected components of the space of framed curves 
in $M$ is identified with the product $\Z_2\times \mathcal C$.
Here the $\Z_2$-factor is the value of the $\spin$-structure on the loop in
the principal $SO(3)$-bundle of $TM$ that corresponds to a framed curve
from the connected component, see~\ref{proofhomotopy}. 
(This value does not depend on the choice of the framed curve in the component.)

The $h$-principle for the Legendrian curves says that the set of the
connected components of the space of Legendrian curves is naturally
identified with the set of homotopy classes of free loops in $CM$ (the
spherical contact bundle of $M$), or which is the same with the set of conjugacy
classes of the elements of $\pi_1(CM)$. Since every contact manifold is
oriented and the contact structure is cooriented, we see that the planes of
the contact structure are naturally oriented. This orientation induces the
orientation of the $S^1$-fibers of $\pr:CM\rightarrow M$. 
Put $f\in\pi_1(CM)$ to be the class of the naturally oriented $S^1$-fiber of 
$\pr:CM\rightarrow M$.

The Theorem of Chow~\cite{Chow} and Rashevskii~\cite{Rashevskii} says
that every connected component of the space of curves contains a Legendrian 
curve. A straightforward verification shows that the insertion 
of the zig-zag into 
the Legendrian curve $K$ (see Figure~\ref{twocusp.fig}) 
changes the value of the $\spin$-structure on the corresponding framed curve. 
It is easy to verify (see~\cite{FuchsTabachnikov}) that the two connected components of the space of Legendrian 
curves that contain $K$ and $K$ with the zigzag correspond respectively to the 
conjugacy classes of $\vec K$ and of $\vec K f$ (or of $\vec K f^{-1}$) in $\pi_1(CM)$. (Here the fact whether
we get $\vec K f$ or $\vec K f^{-1}$ depends on which of the two possible
zig-zags we insert.)

Let $\mathcal L\subset \mathcal F$ be a connected component of the space of
Legendrian curves in $(M, C)$ that corresponds to the conjugacy class of
$\vec K \in \pi_1(CM)$. Then every connected component 
$\mathcal L'\subset \mathcal F$ of the space of Legendrian curves
corresponds to the conjugacy class of $\vec K f^{2n}\in\pi_1(CM)$,
for some $n\in\Z$.

Hence $\mathcal F$ satisfies condition \textrm{II} if and only if for every 
$0\neq n$ the elements $\vec K$ and $\vec K f^{2n}$ are not conjugate in
$\pi_1(CM)$.

Assume that $\mathcal F$ does not satisfy condition \textrm{II}, then there
exists $0\neq n\in\Z$ and $\beta\in\pi_1(CM)$ such that 
\begin{equation}\label{conjugate}
\beta \vec K \beta^{-1}=\vec K f^{2n}\in\pi_1(CM, \vec K(1)).
\end{equation} 

This implies that $\pr_*(\beta)$ and $\pr_*(\vec K)$ commute in $\pi_1(M,
K(1))$. The commutation relation gives a mapping $\mu:T^2\rightarrow M$ of
the two torus $T^2=S^1\times S^1$ such that $\mu |_{(1\times S^1)}=K$ and 
$\mu |_{(S^1\times 1)}=\pr(\beta)$. 

Put $e$ to the Euler class of the contact bundle.
Consider the locally-trivial $S^1$-fibration $p:M'\rightarrow T^2$ induced by 
$\mu$. One can verify that $2n\in \Z=H^2(T^2, \Z)$ is the Euler class of
$p$. On the other hand the Euler class of $p$ is $\mu^*(e)$  and is naturally 
identified with the value of $e$ on the homology class realized by $\mu(T^2)$. 
This implies that if $\mathcal F$ does not satisfy condition \textrm{II}, then there 
exists a homology class $\alpha$ from the statement of the Proposition.

On the other hand the existence of the class $\alpha$ from the statement of
the Proposition implies that there exists a Legendrian curve $K\in\mathcal
F$ such that $\vec K$ is conjugate to $\vec K f^{n}$, for $n$ being the
value of $e$ (the Euler class of the contact bundle) on the homology class 
realized by $\mu(T^2)$. (Proposition~\ref{existcontact} says that 
$e=2\alpha$, for some $\alpha \in H^2(M, \Z)$, and 
hence $n$ is even.) This means that $\mathcal F$ does not satisfy condition
\textrm{II} and we have proved Proposition~\ref{interpretationconditionII}.
\qed

\subsection{Proof of Theorem~\ref{atoroidal}}\label{Proofatoridal}
\begin{emf}\label{proofconditionII}
Similar to the Proof~\ref{ProofinterpretationconditionII} 
of Proposition~\ref{interpretationconditionII} we get that to prove 
that all the components of the space of curves satisfy condition \textrm{II}
it suffices to show that
$\vec K$ and $\vec K f^{2n}$ are not conjugate in $\pi_1(CM)$, for all
$0\neq n\in\Z$ and $\vec K\in\pi_1(CM)$. 
Let $\beta\in\pi_1(CM)$ and
$n\in\Z$ be such that 
\begin{equation}\label{conjugate1}
\beta \vec K \beta^{-1}=\vec K f^{2n}\in\pi_1(CM, \vec K(1)).
\end{equation} 
To prove the statement we have to show that $n=0$.
Identity~\eqref{conjugate1} implies that $\pr_*(\beta)$ and $\pr_* (\vec K)$ commute in $\pi_1(M)$.
Consider a mapping of the two-torus $\mu:T^2=S^1\times S^1\rightarrow M$ 
such that $\mu(S^1\times 1)=\pr (\vec K)$ and $\mu(1\times S^1)=\pr(\beta)$.
(Such a mapping exists since $\pr(\beta)$ and $\pr (\vec K)$ commute in
$\pi_1(M)$.) By the assumption of the Theorem
$\mu:\pi_1(T^2)\rightarrow\pi_1(M)$ has a nontrivial kernel. Thus there
exist $i,j\in\Z$ with at least one of $i$ or $j$ being nonzero such that 
$\pr(\vec K)^i=\pr(\beta)^j \in \pi_1(M)$, and hence
\begin{equation}
\vec K^i=\beta ^j f^l, \text{ for some } l\in \Z. 
\end{equation}

Since the situation is symmetric, we assume that $j\neq 0$.

Thus $\vec K^i \vec K^i= \vec K^i\beta ^j f^l=\beta ^j f^l\vec K^i$.
Applying~\eqref{conjugate1} to the last identity we get that $f^{2nij}=1$. 
Since $\pi_2(M)=0$ we see that $f$ has infinite order in $\pi_1(CM)$, and
hence $2nij=0$. If $n$ is zero, then we are done. Hence we have to look at
the case of $i=0$. (We assumed that $j\neq 0$.) 

If $i=0$, then $\pr_*(\beta ^j)=1$. Since we assumed that $j\neq 0$ the
assumption of the Theorem implies that $\pr_*(\beta)=1$. Thus $\beta$ is a
power of $f$ and hence is in the center of $\pi_1(CM)$. Thus $n=0$ in and we have 
proved that all the components of the space of framed curves satisfy condition \textrm{II}.
\end{emf}

\begin{emf}\label{ProofconditionI} 
The proof of the fact that every component $\mathcal F$ of the space 
of framed curves satisfies condition \textrm{I} is obtained by a straightforward 
generalization of the proof of Theorem~\ref{fibration}
(see~\ref{Prooffibration}) and in fact is much simpler. 

From the condition that for every $\mu:T^2\rightarrow M$ the homomorphism
$\mu_*:\pi_1(T^2)\rightarrow \pi_1(M)$ is not injective one gets that 
for every $\beta\in Z(K_1)<\pi_1(M, K_1(1))$ there exist $n,i\in\Z$ such that 
at least one of $n$ and $i$ is nonzero and $\beta^n=K_1^i\in\pi_1(M, K_1(1))$.

If the connected component of the space of unframed curves that contains
$K_1$ consists of noncontractible loops, then $n\neq 0$ (otherwise 
$K_1\in\pi_1(M, K_1(1))$ is an element of finite order which contradicts to the
assumptions of the Theorem). And then a certain nonzero degree of every
$\beta\in Z(K_1)$ can be realized as $t_1(\gamma_2^i)=t_2(\gamma_2^i)$, for
the loop $\gamma_2$ introduced in~\ref{notcontractible}. Similar
to~\ref{notcontractible} this finishes the proof in the case of curves from
$\mathcal F$ realizing noncontractible loops.

The proof in the case of $\mathcal C$ consisting of contractible loops is also 
completely clear. This finishes the proof of Theorem~\ref{atoroidal}.\qed
\end{emf}

\subsection{Proof of Theorem~\ref{example1}}\label{Proofexample1}
\begin{emf}{\em Proof of statement {\bf a:} of Theorem~\ref{example1}.\/}
Clearly (see Figure~\ref{homotopy.fig}) the two Legendrian knots $K_0$ and
$K_1$ belong to the same component of the space of Legendrian curves. It is easy to see 
that if $K_0$ realizes the isotopy class of the framed knot $\tilde K_0$, then
$K_1$ realizes the isotopy class of $\tilde K_0^{-2}$ (see~\ref{obstruction} for
the definition of $\tilde K_0^{-2}$).
Below we show that $\tilde K_0$ and $\tilde K_0^{-2}$ are isotopic framed
knots. 

Let $t\times S^2$ be the sphere that crosses $\tilde K_0$ at exactly one
point, and let $N=[0,1]\times S^2$ be the tubular neighborhood of $t\times
S^2$. Fix $x\in S^2$ (below called the North pole) and the direction in $T_x
S^2$ (below called the zero meridian). 
We can assume that the knot $\tilde K_0$ inside $N=[0,1]\times S^2$ looks as
follows: it intersects each $y\times S^2\subset [0,1]\times S^2$ at the
North pole of the corresponding sphere, and the framing of the knot is parallel to the zero meridian.

Consider an automorphism $\nu:S^1 \times  S^2\rightarrow S^1 \times  S^2$
that is identical outside of $N=[0,1]\times S^2$ such that it rotates each 
$y\times S^2\in [0,1]\times S^2$ by $4\pi y$ around the North pole in the clockwise 
direction. Clearly under this automorphism $\tilde K_0$ gets two extra negative 
twists of the framing and $\nu(\tilde K_0)=\tilde K_0^{-2}$. On the other hand 
it is easy to see that $\nu$ is diffeotopic to the identity, since it corresponds 
to the contractible loop in $SO(3)=\R P^3$. Hence we see that $\tilde K_0$ and 
$\tilde K_0^{-2}$ are isotopic framed knots. This finishes the proof of the 
statement {\bf a:} of Theorem~\ref{example1}.
\end{emf}

To prove statement {\bf b:} of the Theorem we need the following
Proposition.

\begin{prop}\label{propS1xS2}
Let $C$ be a cooriented contact structure on $M=S^1\times S^2$
with a nonzero Euler class $e$ of the contact bundle. Let $CM$ be the spherical contact
bundle, let $\pr:CM\rightarrow M$ be the corresponding locally trivial 
$S^1$-fibration, and let $f\in\pi_1(CM)$ be the class of an oriented $S^1$-fiber 
of $\pr$. Then $f$ is of finite order in $\pi_1(CM)$ and $\pi_2(CM)=0$.
\end{prop}

\begin{emf}{\em Proof of Proposition~\ref{propS1xS2}.\/}
Consider the oriented $2$-plane bundle $p:\xi\rightarrow S^2$ that is the
restriction of the contact bundle over $M$ to a sphere 
$1\times S^2\subset S^1\times S^2$. The Euler class of $\xi$ is the value of 
$e$ on the homology class realized by $1\times S^2$, and hence is nonzero.
Let $S\xi$ be the manifold obtained by the fiberwise spherization of $p$, and
let $\bar p:S\xi\rightarrow S^2$ be the corresponding locally trivial
$S^1$-fibration. Since the Euler class of $\xi$ is nonzero we get that 
a certain multiple of the class of the fiber of $\bar p$ is homologous to
zero. But $\pi_1(S\xi)$ is generated by the class of the fiber 
and hence the class of the fiber of $\bar p$ is of finite order in
$\pi_1(S\xi)$. This implies that $f\in\pi_1(CM)$ is of finite order.

The statement that $\pi_2(CM)=0$ follows from the exact homotopy sequence of 
$\pr:CM\rightarrow M$ and the fact that $f\in\pi_1(CM)$ is of finite order.
\qed
\end{emf}  

\begin{emf}{\em Proof of statement {\bf b:} of Theorem~\ref{example1}.\/}
Let $\mathcal L$ be the connected component of the space of Legendrian
curves that contains $K_0$ and $K_1$.
Figure~\ref{homotopy.fig} shows that $K_0$ can be changed to $K_1$ (in the
space of Legendrian curves) by a
sequence of isotopies and one passage through a transverse double point. Hence if there 
exists a $\Z$-valued invariant $I$ of Legendrian knots from $\mathcal L$
that increases by one under every positive passage through a transverse double point of a Legendrian 
knot, then it distinguishes $K_0$ and $K_1$. (Clearly if such $I$ does
exist, then it is an order one invariant of Legendrian knots.) Below we show the existence of
such $I$ in the connected component $\mathcal L$ 
of the space of Legendrian curves that contains $K_0$.

\begin{figure}[htbp]
 \begin{center}
  \epsfxsize 10cm
  \hepsffile{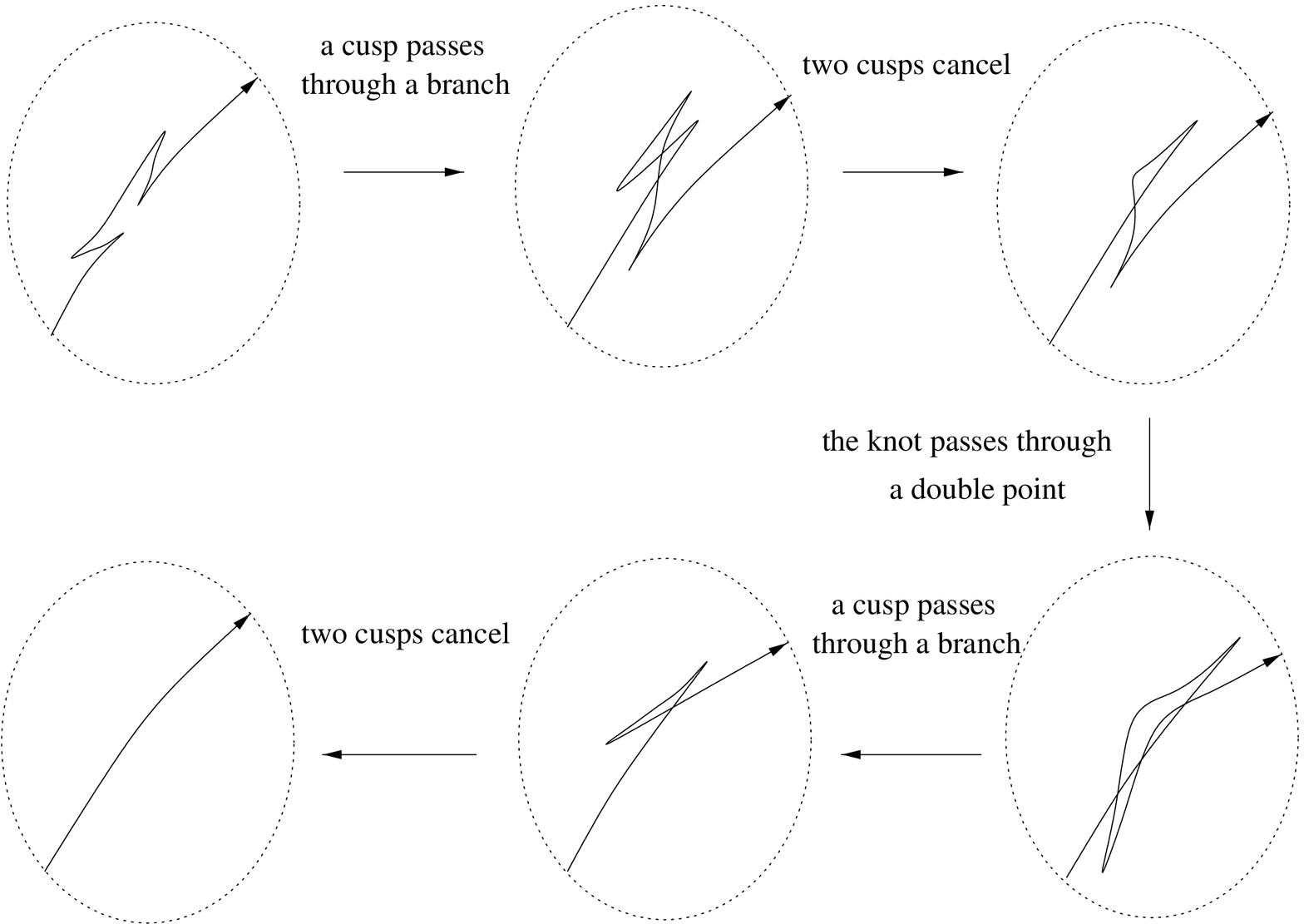}
 \end{center}
\caption{}\label{homotopy.fig}
\end{figure}

Put $I(K_0)=0$. Let $K'\in\mathcal L$ be a Legendrian knot, and let 
$\gamma$ be a generic path in $\mathcal L$ connecting $K_0$ and $K'$. Let
$J_{\gamma}$ be the set of moments when $\gamma$ crosses the discriminant (i.e.
a subspace of singular knots) in $\mathcal L$, and let $\sigma_j$, $j\in
J_{\gamma}$, be the signs of these crossings. For a generic path $\gamma\subset \mathcal L$ put 
$\Delta_I(\gamma)=\sum_{j\in J_{\gamma}}\sigma_j$. It is clear that if $I$ (with
$I(K_0)=0$) does exist,
then $I(K')=\Delta_I(\gamma)$. To show that $I$ does exist we have to
verify that $\Delta_I(\gamma)$ does not depend on the choice of a generic
path $\gamma$ connecting $K_0$ and $K'$, or which is the same we have to show
that $\Delta_I(\alpha)$ is zero for every generic closed loop $\alpha$ 
connecting  $K_0$ to itself. 

There are two codimension two strata of the discriminant of $\mathcal L$.
They are formed respectively by singular Legendrian knots with two 
transverse double points, and by Legendrian knots with one double point at
which the two intersecting branches are tangent of order one.
It is easy to see that $\Delta_I(\beta)=0$, for every small closed loop
$\beta$ going around a codimension two stratum of $\mathcal L$. 

This implies that $\Delta_I(\gamma)$ depends only on the element of
$\pi_1(\mathcal L, K_0)$ realized by $\gamma$. Hence to prove the
existence of $I$ it suffices to show that $\Delta_I(\gamma)=0$ for every
$\gamma\in\pi_1(\mathcal L, K_0)$.

Clearly $\Delta_I(\gamma^p)=p\Delta_I(\gamma)$ and since $\Z$ is torsion
free, we get that to prove Theorem~\ref{example1} it suffices to show that 
$\Delta_I(\gamma^p)=0$ for a certain nonzero $p\in \Z$.

The $h$-principle says that the space of Legendrian curves in $(M,C)$ is weak homotopy equivalent to 
the space $\Omega CM$ of free loops in the spherical contact bundle $CM$ of 
$M=S^1\times S^2$. (The mapping giving the equivalence lifts a Legendrian
curve $K$ to a loop $\vec K$ in $CM$ by sending $t\in S^1$ to the point in
$CM$ that corresponds to the direction of the velocity vector of $K$ at $K(t)$.)

Hence $\pi_1(\mathcal L, K_0)$ is naturally isomorphic to $\pi_1(\Omega CM,
\vec K_0)$. Since $\pi_2(CM)=0$ by Proposition~\ref{propS1xS2}, we get (see~\ref{Hansen}) that $\pi_1(\Omega
CM, \vec K_0)$ is isomorphic to $Z(\vec K_0)$ the centralizer of $\vec
K_0\in \pi_1(CM, \vec K_0(1))$. Using Propositions~\ref{commute}
and~\ref{propS1xS2} we see that $\pi_1(CM)=\Z\oplus\Z_p$, for some nonzero $p\in\N$, 
or $\pi_1(CM)=\Z$. Hence there exists $n\in\Z$ and nonzero $m\in\Z$ 
such that $\gamma^m=\vec K_0^n\in\pi_1(CM, \vec K_0(1))$. (One should take
$n$ and $m$ to be divisible by $p$ if $\pi_1(CM)=\Z\oplus\Z_p$.)
But the loop
$\alpha$ in $\pi_1(\mathcal L, K_0)$ corresponding to $\vec K_0^n$ is just the sliding
$n$-times of $K_0$ along itself according to the orientation.
(This deformation is  induced by the rotation 
of the parameterizing circle.) This loop does not intersect the
discriminant, 
and hence $\Delta_I(\alpha)=0$. This finishes the proof of 
the statement {\bf b:} of Theorem~\ref{example1}.\qed
\end{emf}

\subsection{Proof of Theorem~\ref{example2}}\label{Proofexample2}
\begin{emf}\label{part1}{\em $K_1$ and $K_2$ are homotopic Legendrian
curves, and $K_1$ and $K_2$ realize isotopic framed knots.\/}
Let $f_1\in\pi_1(CM)$ be the class of the $S^1$-fiber of the fibration 
$\pr:CM\rightarrow M$. The $h$-principle says that 
the connected component of the space of Legendrian curves that contains $K$
corresponds to the conjugacy class of $\vec K\in\pi_1(CM)$. From the work of Fuchs and
Tabachnikov it is easy too see that the connected components containing 
$K_1$ and $K_2$ correspond to the conjugacy classes of $\vec K_1=\vec K
f_1^r$ 
and  of $\vec K_2=\vec K f_1^{-r}$. 
Let $f_2\in\pi_1(M)$ be an element projecting to the class $f\in\pi_1(M)$ of 
the $S^1$-fiber of $p:M\rightarrow F$. The value of the Euler class of the 
contact bundle  on the homology class realized by $\mu(T^2)$ is equal to 
$2r\in\Z$. (Here $\mu$ is the mapping from the description of the Euler
class of the contact bundle.) And because of the reasons explained in the proof of
Proposition~\ref{interpretationconditionII},
we get that $\vec K f_2=f_2 \vec K f_1^{2r}$. Now Proposition~\ref{commute}
implies that $\vec K_1$ and $\vec K_2$ are conjugate in $\pi_1(CM)$ and
hence $K_1$ and $K_2$ are in the same component of the space of Legendrian
curves. 

The fact that $K_1$ and $K_2$ realize isotopic framed knots is clear,
because as unframed knots they are the same, and as it is shown
in~\cite{FuchsTabachnikov} every pair of extra cusps corresponds to the
negative extra turn of the framing.
\end{emf}

\begin{emf}\label{ideaexample2}{\em Idea of the proof of of the fact that $K_1$ and
$K_2$ can be distinguished by an order one invariant of Legendrian knots.\/}  
Let $d$ be a point in $M$. Let $K_s$ be a singular unframed knot with one 
double point.  The double point separates $K_s$ into two oriented loops. 
Deform $K_s$ preserving the double point, so that the double point is located
at $d$. Choosing one of the two loops of $K_s$ we obtain
an ordered set of two elements $\delta_1, \delta_2\in \pi_1(M,d)$, or which
is the same an element $\delta_1\oplus\delta_2\in\pi_1(M,
d)\oplus\pi_1(M,d)$. 
Clearly there is a unique element of the set $B$ that corresponds to the
original singular unframed knot $K_s$, where $B$ is is the quotient set of 
$\pi_1(M,d)\oplus\pi_1(M,d)$ via the consequent actions of the following groups:
\begin{description}
\item[1] $\pi_1(M)$ whose element $\xi$ acts on
$\delta_1\oplus\delta_2\in\pi_1(M)\oplus\pi_1(M)$ by sending it to 
$\xi^{-1} \delta_1 \xi\oplus \xi^{-1} \delta_2
\xi\in\pi_1(M)\oplus\pi_1(M)$.
(This corresponds to the ambiguity in deforming $K_s$, so
that the double point is located at $d$.) 

\item[2] $\Z_2$ that acts via the cyclic permutation of the two summands.
(This corresponds to the ambiguity in the choice of one of the two
loops of $K_s$.)
\end{description}

Thus we have a mapping $\nu$ from the set of singular unframed knots with 
one double point to $B$. Let $\alpha:B\rightarrow \Z$ be the function such
that 
\begin{description}
\item[a] $\alpha(b)=0$, provided that $b$ contains the class of
$1\oplus\delta\in\pi_1(M)\oplus\pi_1(M)$, for some $\delta\in\pi_1(M)$,
\item[b] $\alpha(b)=1$ otherwise.
\end{description}

Assume that $I^{\mathcal L}$ is an invariant of Legendrian knots 
from $\mathcal L$ such that under every (generic transverse) 
positive passage through a
discriminant in $\mathcal L$ it increases by $\alpha\circ\nu (K_s)$, where $K_s$ is the
unframed singular knot corresponding to the crossing of the discriminant.
Clearly such 
$I^{\mathcal L}$ is an order one invariant of framed knots from 
$\mathcal L$. To prove the Theorem we show the existence of such
$I^{\mathcal L}$, and then we show that it distinguishes 
$K_1$ and $K_2$.
\end{emf}

\begin{emf}\label{existence}{\em The existence of $I^{\mathcal L}$.\/}
Let $\gamma$ be a generic closed loop in $\mathcal L$ that starts with $K_1$.
Let $J_{\gamma}$ be the set of instances when $\gamma$ crosses the discriminant (i.e.
a subspace of singular knots) in $\mathcal L$ and let $\sigma_j$, $j\in
J_{\gamma}$, be the
signs of these crossings. Let $J'_{\gamma}\subset J_{\gamma}$ be those instances for 
which the value of $\alpha\circ \nu$ on the corresponding singular unframed knots is $1$.
For a generic path $\gamma\subset \mathcal L$ put 
$\Delta^{\mathcal L}_I(\gamma)=\sum_{j'\in J'_{\gamma}}\sigma_{j'}$. 

Let $\mathcal C$ be the connected component of the space of unframed curves 
obtained by forgetting the framings on curves from $\mathcal F$, and let
$K'_1$ be the unframed knot obtained by forgetting the framing on $K_1\in
\mathcal L\subset \mathcal F$. Similarly to the above for a generic closed 
loop $\gamma$ in $\mathcal C$ starting with $K'_1$ we put 
$\Delta^{\mathcal C}_I(\gamma)=\sum_{j'\in J'_{\gamma}}\sigma_{j'}$. (As above 
$J'_{\gamma}$ is the set of instances when the value of $\alpha\circ \nu$ on
the singular unframed knots obtained under $\gamma$ is equal to $1$, and $\sigma_{j'}$,
$j'\in J'_{\gamma}$ are the signs of the corresponding crossings of the discriminant.)

Similarly to~\ref{Proofexample1} we get that to prove the existence of $I^{\mathcal L}$
it suffices to show that $\Delta_I^{\mathcal L}(\gamma)=0$, for every generic closed loop $\gamma$. 

Similarly to~\ref{Proofexample1} and to~\ref{mainidea}
we get that $\Delta_I^{\mathcal L}:\pi_1(\mathcal L,
K_1)\rightarrow \Z$ and $\Delta_I^{\mathcal C}:\pi_1(\mathcal C,
K'_1)\rightarrow \Z$ are homomorphisms. Moreover it is clear that if 
$\gamma'\in\pi_1(\mathcal C)$ is the element corresponding to 
$\gamma\in\pi_1(\mathcal L)$, then we have 

\begin{equation}\label{gammagamma'}\Delta_I^{\mathcal L}(\gamma)=
\Delta_I^{\mathcal C}(\gamma').
\end{equation}

The $h$-principle says that the space of Legendrian curves in $(M,C)$ is
weak homotopy equivalent to the space of free loops in $CM$ the spherical
contact bundle of $M$. (The mapping that gives an equivalence lifts a
Legendrian curve $K$ in $(M, C)$ to a loop $\vec K$ in $CM$ by mapping $t\in
S^1$ to the point of $CM$ that corresponds to the velocity vector of $K$ at
$K(t)$.) Since $\pi_2(CM)=0$ for $M$ from the statement of the Theorem, 
we obtain (see~\ref{Hansen}) the natural isomorphism 
$t: \pi_1(\mathcal L, K_1)\rightarrow Z(\vec K_1)< \pi_1(CM, \vec K_1(1))$.
Since $\Delta^{\mathcal L}_I(\gamma^p)=p\Delta^{\mathcal L}_I(\gamma)$ and $\Z$ is torsion
free, we get that to show the existence of $I^{\mathcal L}$ it suffices to show 
that for every $\beta\in Z(\vec K_1)<\pi_1(CM,\vec K_1(1))$ there exist $0\neq n\in\Z$ 
and $\gamma\in\pi_1(\mathcal L, K_1)$ such that $t(\gamma)=\beta^n\in\pi_1(CM, \vec K_1(1))$ and
$\Delta^{\mathcal L}_I(\gamma)=0$.

Let $f\in\pi_1(M, K_1(1))$ be the class of the $S^1$-fiber of $p:M\rightarrow
F$. Let $f_1\in\pi_1(CM)$ be the class of an oriented $S^1$-fiber of
$pr:CM\rightarrow M$, and let $f_2$ be an element of $\pi_1(CM)$ such that 
$pr_*(f_2)=f\in\pi_1(M, K_1(1))$. 

Take $\beta\in Z(\vec K_1)$, then $\pr_*(\beta)\in Z(K_1)$.
Proposition~\ref{toughandtechnical} implies that there exist $0\neq n\in\Z$
and $i,j\in\Z$ such that $K_1^if^j=(\pr_*(\beta))^n\in\pi_1(M, K_1(1))$.
Using Proposition~\ref{commute} we get that 
\begin{equation}\label{expression}
\beta^n=\vec K_1^if_2^j f_1^l\text{ for some }i,j,l\in\Z. 
\end{equation}

As it was explained in~\ref{part1} we have 
\begin{equation}\label{extraobstruction}
\vec K_1 f_2= f_2 \vec K_1 f_1^{2r}.
\end{equation} 

Since $\beta\in Z(\vec K_1)$ we get that $\vec K_1\beta^n=
\beta^n \vec K_1$, and using~\eqref{expression} we see that 
$\vec K_1 \vec K_1^if_2^j f_1^l=\vec K_1^if_2^j f_1^l\vec K_1$.
Using~\eqref{extraobstruction}, Proposition~\ref{commute}, 
and the fact that $f_1$ has infinite order in $\pi_1(CM)$, we see that $j=0$
in~\eqref{expression}.

Hence 
\begin{equation}\label{newexpression}
\beta^n=\vec K_1^i f_1^l,\text{ for some }i,l\in\Z. 
\end{equation}

Clearly $f_1\in Z(\vec K_1)$ and hence by the $h$-principle 
there is a loop
$\gamma_3\in\pi_1(\mathcal L, K_1)$ such that $t(\gamma_3)=f_1$. 
Let $\gamma_2\in\pi_1(\mathcal L, K_1)$ be the loop corresponding to the
deformation under which $K_1$ slides once around itself according to the
orientation of $K_1$. (This deformation is induced by the rotation of the
circle parameterizing $K_1$.) Clearly $\Delta^{\mathcal L}_{I}(\gamma_2)=0$.

Thus to prove the existence of $I^{\mathcal L}$ 
it suffices to show that $\Delta^{\mathcal L}_{I}(\gamma)=0$, for
$\gamma\in\pi_1(\mathcal L, K_1)$ such that $t(\gamma)=\vec K_1^i f_1^l$.
But this $\gamma$ is $\gamma_2^i \gamma_3^l$. Thus to prove the Theorem it
suffices to show that $0=\Delta^{\mathcal L}_{I}(\gamma)=
i\Delta^{\mathcal L}_{I}(\gamma_2)+l\Delta^{\mathcal L}_{I}(\gamma_3)$. Since 
$\Delta^{\mathcal L}_{I}(\gamma_2)=0$
we get that $\Delta^{\mathcal L}_{I}(\gamma)=l\Delta^{\mathcal L}_{I}(\gamma_3)$ and thus it suffices
to show that $\Delta^{\mathcal L}_{I}(\gamma_3)=0$.

The loop $\gamma_3\in\pi_1(\mathcal L, K_1)$ corresponds to a loop 
$\gamma'_3\in\pi_1(\mathcal C, K'_1)$. Identity~\eqref{gammagamma'} says that 
$\Delta^{\mathcal L}_{I}(\gamma_3)=\Delta^{\mathcal C}_{I}(\gamma'_3)$.
Hence we have to show that $\Delta^{\mathcal C}_{I}(\gamma'_3)=0$.

Using the fact that $\pr_*(f_1)=1\in\pi_1(M, K_1(1))$, 
the $h$-principle for curves~\ref{h-principleforcurves}, 
and observations of~\ref{mainidea}, and of~\ref{obstructiondecrease} we obtain
that $\gamma'_3\in\pi_1(\mathcal C, K'_1)$ can be realized by a sequence
of deformations shown in Figure~\ref{obstruction.fig}. One of the loops of a
singular knot arising under the deformation shown in
Figure~\ref{obstruction.fig} is contractible and the value of $\alpha$ on such
a singular knot was put to be zero. Thus we have that $\Delta^{\mathcal
C}_{I}(\gamma'_3)=0$, and this finishes the Proof of the existence of
$I^{\mathcal L}$.
\end{emf}

\begin{emf}\label{distinguish}
{\em Let us show that $I^{\mathcal L}$ distinguishes $K_1$ and $K_2$.\/}
As before let $K'_1$ be the unframed knot obtained by forgetting the framing
on $K_1$. Let $\rho:S^1\rightarrow \mathcal L$ be a generic path connecting 
$K_1$ and $K_2$. To prove the Theorem we have to show that 
$\Delta_{I^{\mathcal L}}(\rho)\neq 0$.

The homotopy $\rho$ gives rise to a homotopy $\bar \rho:S^1\rightarrow \mathcal C$ 
connecting $K'_1$ to itself. This $\bar \rho$ has the property that the value
of the Euler class of the contact bundle on the homology class realized by
the corresponding mapping $T^2\rightarrow M$ is nonzero. Moreover it is clear that 
$\Delta_{I^{\mathcal C}}(\bar \rho)=\Delta_{I^{\mathcal L}}(\rho)$. 

Using the usual arguments we get that to prove the Theorem it suffices to 
show that for every $\gamma\in\pi_1(\mathcal C, K'_1)$ there exists 
$n\neq 0$ such that either 
\begin{description}
\item[a] the value of the Euler class of the contact bundle on the homology
class realized by the mapping of $T^2$ corresponding to $\gamma^n$ is zero,
or
\item[b] $\Delta_{I^{\mathcal C}}(\gamma^n)\neq 0$.
\end{description}
It is easy to see that the element of $\pi_1(\mathcal C, K'_1)$ realized by the 
deformation in Figure~\ref{obstruction.fig} is in the center of 
$\pi_1(\mathcal C, K'_1)$. Now from the proof of Theorem~\ref{fibration} one concludes 
that for every $\gamma\in\pi_1(\mathcal C, K'_1)$ there exist $0\neq n\in\Z$ such 
that $\gamma^n$ can be expressed as a product $\gamma^n=\gamma_1^i \gamma_2^j
\gamma_4^k$ of the powers of the loops $\gamma_1$, $\gamma_2$ 
described in~\ref{notcontractible} and of the loop
$\gamma_4$ that is the deformation described in Figure~\ref{obstruction.fig}.

It is easy to see that if $i\neq 0$ then $\Delta_{I^{\mathcal C}}\neq
0$, and that if $i=0$ then the value of the Euler class on the homology
class realized by the corresponding mapping of $T^2$ is zero. This finishes
the proof of Theorem~\ref{example2}.
\qed
\end{emf}

{\bf Acknowledgements.} 
I am very grateful to O.~Ya.~Viro for the valuable discussions and
suggestions. I am deeply thankful to H.~Geiges, A.~Stoimenow and
S.~Tabachnikov for the valuable suggestions, and to O.~Baues, M.~Bhupal, 
N.~A'~Campo, A.~Cattaneo, J.~Fr\"ohlich, A.~Shumakovich and V.~Turaev
for many valuable discussions.

This paper was written during my stay at the ETH Zurich and I would like to
thank the staff of the ETH for providing the excellent working conditions.

\end{document}